\documentclass[12pt]{iopart}

\expandafter\let\csname equation*\endcsname\relax
\expandafter\let\csname endequation*\endcsname\relax
\usepackage{amsmath}
\UseRawInputEncoding
\usepackage{graphicx}
\usepackage{subfigure}
\usepackage[ruled, linesnumbered]{algorithm2e}
\usepackage{setspace}
\usepackage{array}
\usepackage{booktabs}
\usepackage{tabularx}
\usepackage{multirow}
\usepackage{hyperref}
\usepackage{cleveref}
\hypersetup{hypertex=true,
	colorlinks=true,
	linkcolor=blue,
	anchorcolor=blue,
	citecolor=blue}
\usepackage{graphbox}

\usepackage{appendix}
\newcommand*{\fullref}[1]{\nameref{#1}}

\begin{document}

\title[]{Fast Iterative Reconstruction for Multi-spectral CT by a Schmidt Orthogonal Modification Algorithm (SOMA) }

\author{Huiying Pan$^1$, Shusen Zhao$^2$, Weibin Zhang$^1$, Huitao Zhang$^{1,2}$ and Xing Zhao$^{1,2,*}$}

\address{$^1$ School of Mathematical Sciences, Capital Normal University, Beijing, 100048, China}
\address{$^2$ National Center for Applied Mathematics Shenzhen(NCAMS), Southern University of Science and Technology, Shenzhen, 518055, China}
\ead{* zhaoxing$\_$1999@126.com}

\vspace{10pt}
\begin{indented}
\item[Received] xxxxxx
\item[Accepted for publication] xxxxxx
\item[Published] xxxxxx
\end{indented}

\begin{abstract}
Multi-spectral CT (MSCT) is increasingly used in industrial non-destructive testing and medical diagnosis because of its outstanding performance like material distinguishability. The process of obtaining MSCT data can be modeled as nonlinear equations and the basis material decomposition comes down to the inverse problem of the nonlinear equations. For different spectra data, geometric inconsistent parameters cause geometrical inconsistent rays, which will lead to mismatched nonlinear equations. How to solve the mismatched nonlinear equations accurately and quickly is a hot issue. This paper proposes a general iterative method to invert the mismatched nonlinear equations and develops Schmidt orthogonalization to accelerate convergence. The validity of the proposed method is verified by MSCT basis material decomposition experiments. The results show that the proposed method can decompose the basis material images accurately and improve the convergence speed greatly.

\vspace{1pc}
\end{abstract}

\section{Introduction}
\label{sect1}
Computed Tomography (CT) can show the internal details without destroying or damaging the objects and has been widely used in many fields such as medicine \cite{ref1,ref2,ref3}, materials \cite{ref4,ref5}, geological engineering \cite{ref6} and so on \cite{ref7,ref8}. Multi-spectral CT (MSCT) takes photon energy into account \cite{ref9} and obtains more information about the objects \cite{ref10}. Compared with conventional CT, MSCT has better artifact removing performances \cite{ref11,ref12}, quantitative detectability \cite{ref13,ref14} and material distinguishability \cite{ref15,ref16,ref17}. So MSCT is increasingly used in industry \cite{ref18} and medicine \cite{ref19}, especially in medical diagnosis \cite{ref20,ref21,ref22}.

\begin{figure}[htbp]
	\centering
	\subfigure[]{
		\includegraphics[width=4.5cm]{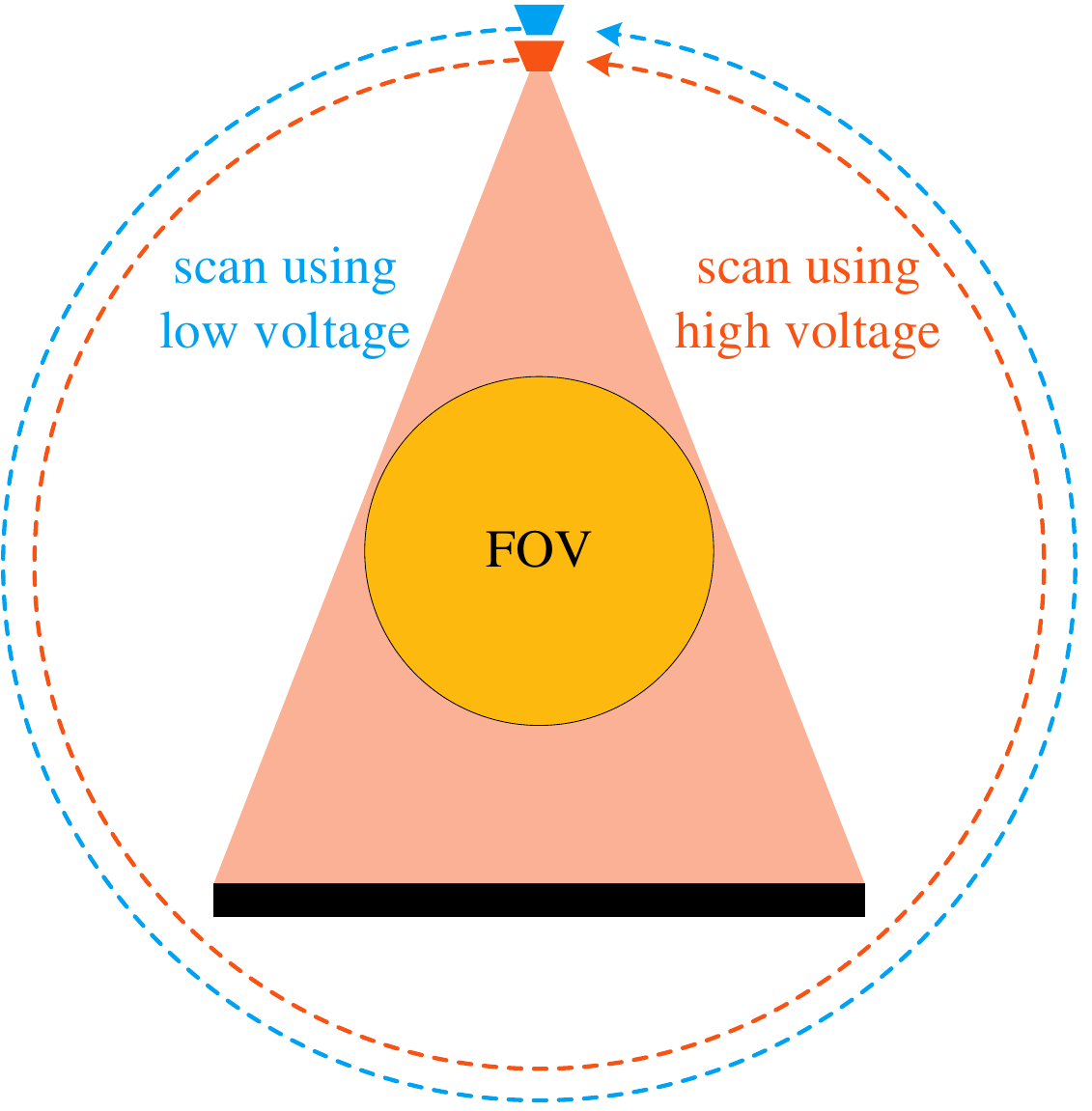}
		\label{fig1_a}}
	\subfigure[]{
		\includegraphics[width=4.5cm]{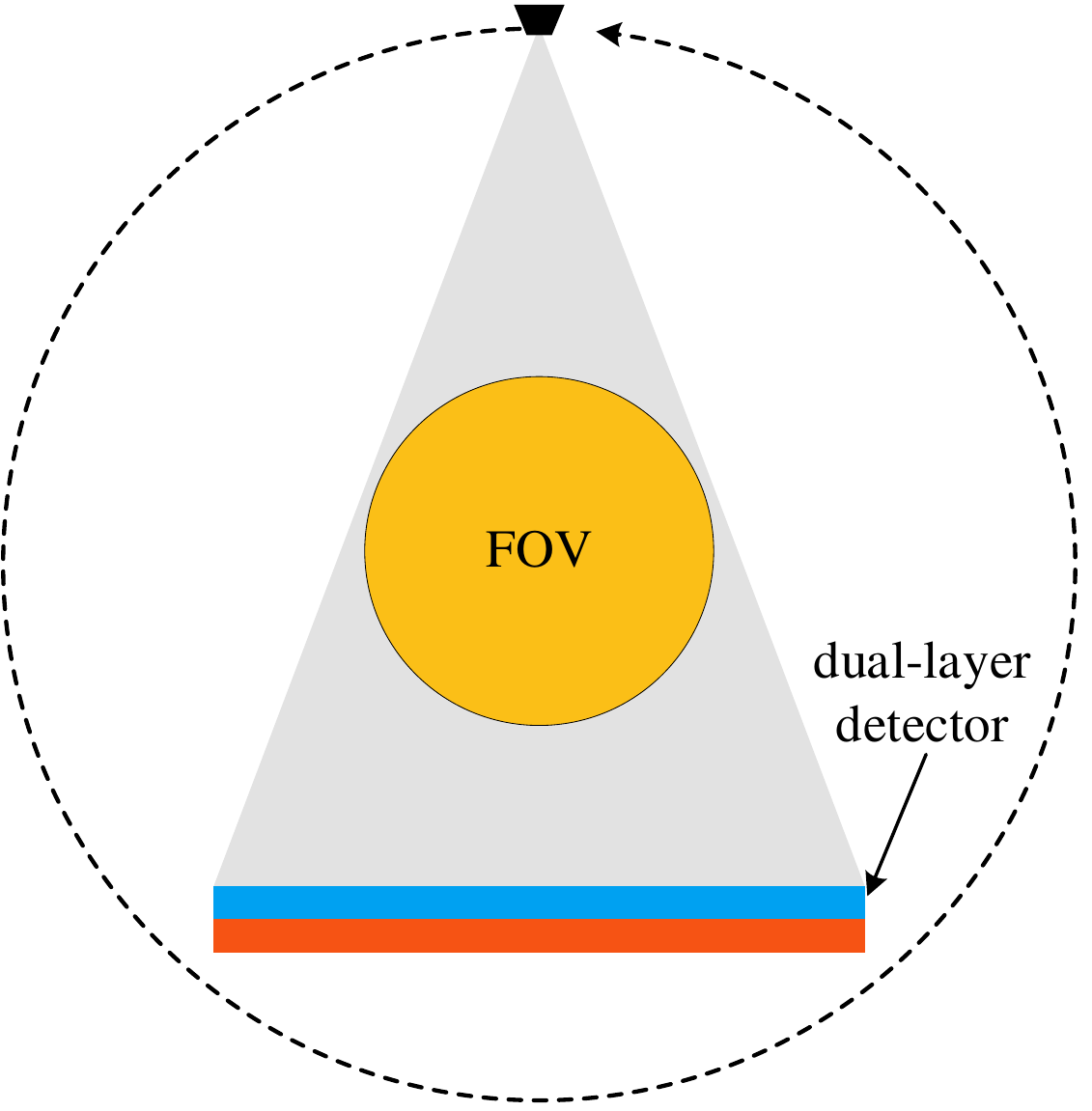}
		\label{fig1_b}}
	\subfigure[]{
		\includegraphics[width=4.5cm]{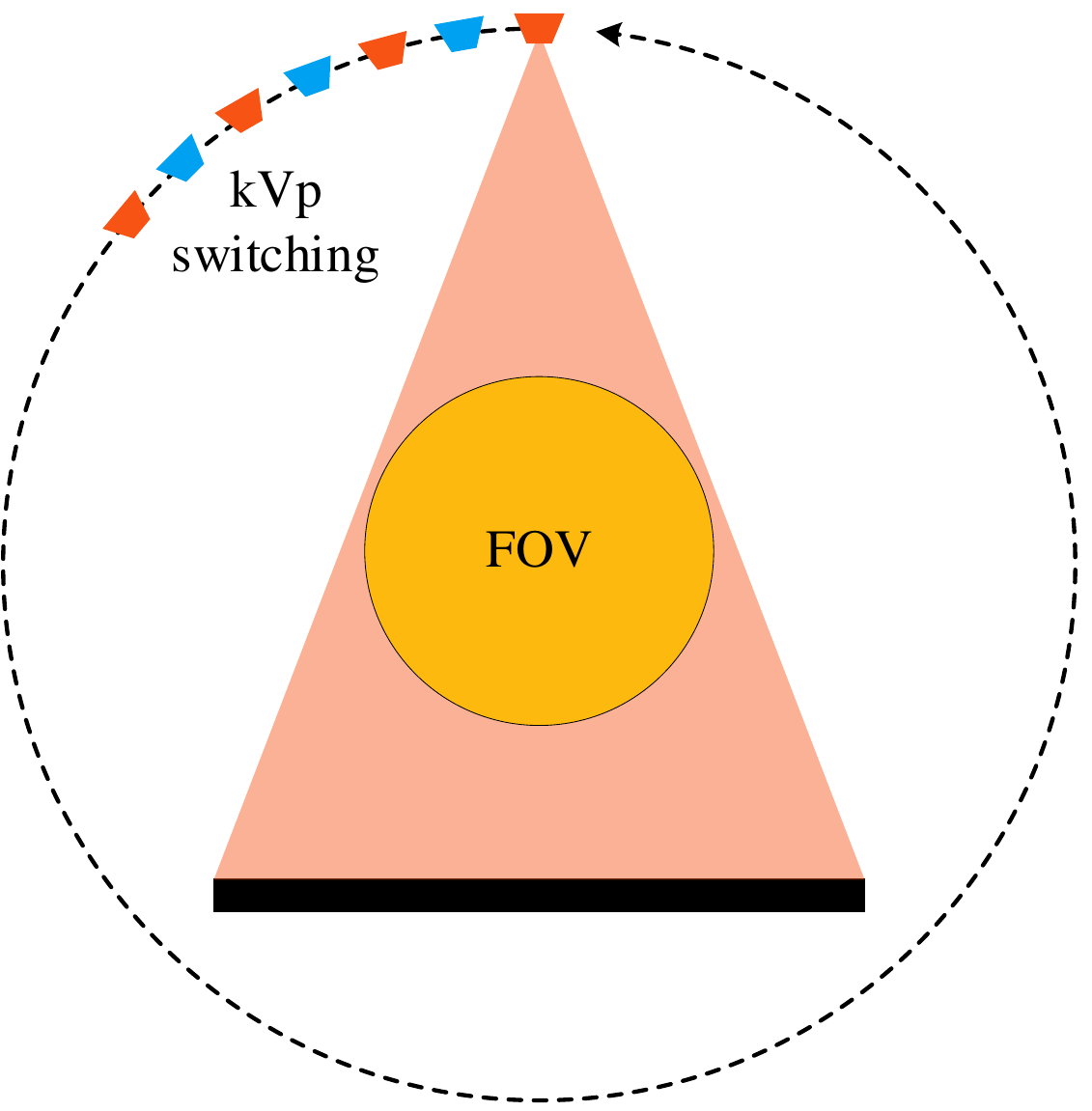}
		\label{fig1_c}} \hspace{-2cm}
	\vspace{-2.5mm}
	
	\subfigure[]{
		\includegraphics[width=4.5cm]{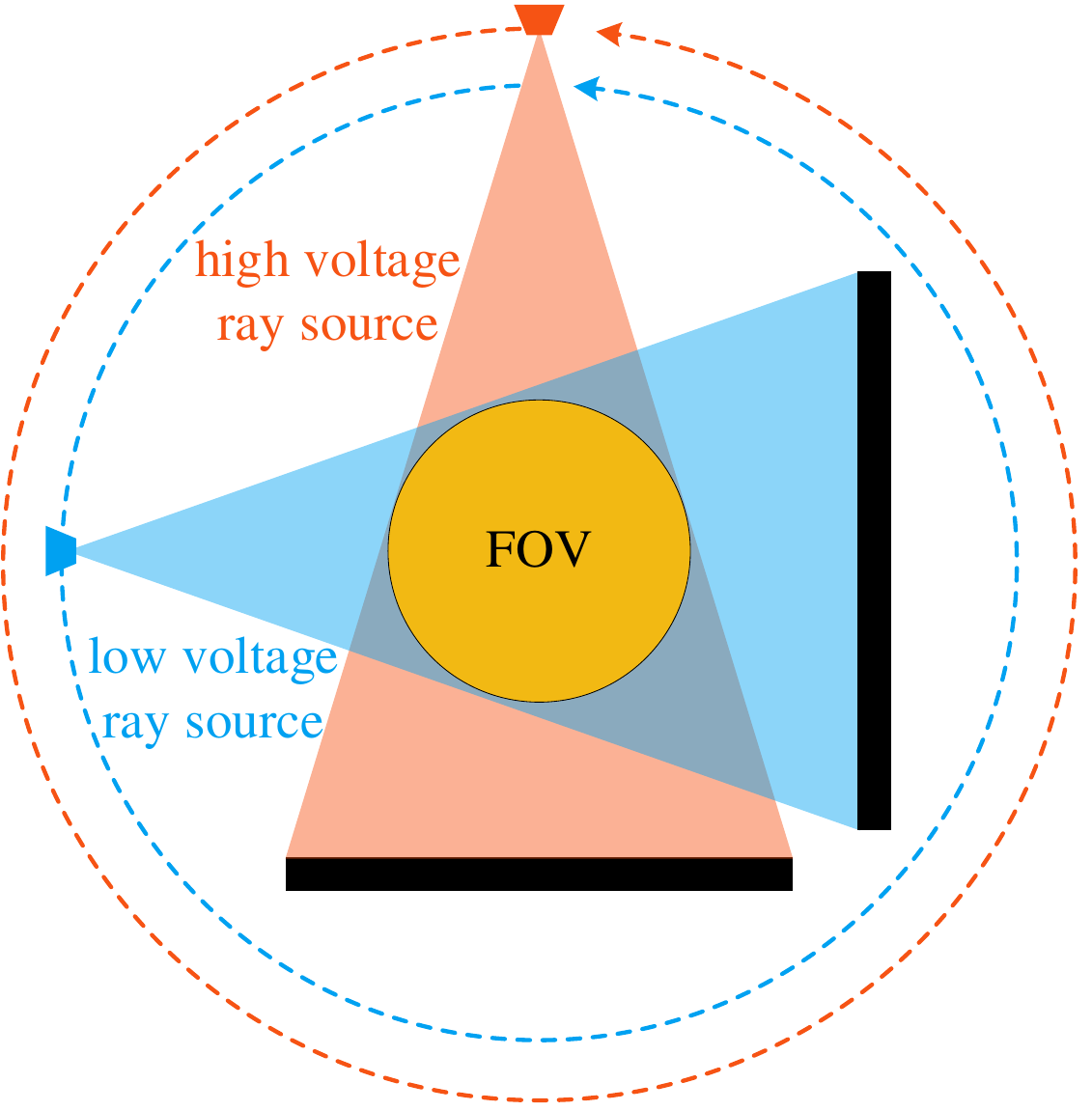}
		\label{fig1_d}}
	\subfigure[]{
		\includegraphics[width=4.5cm]{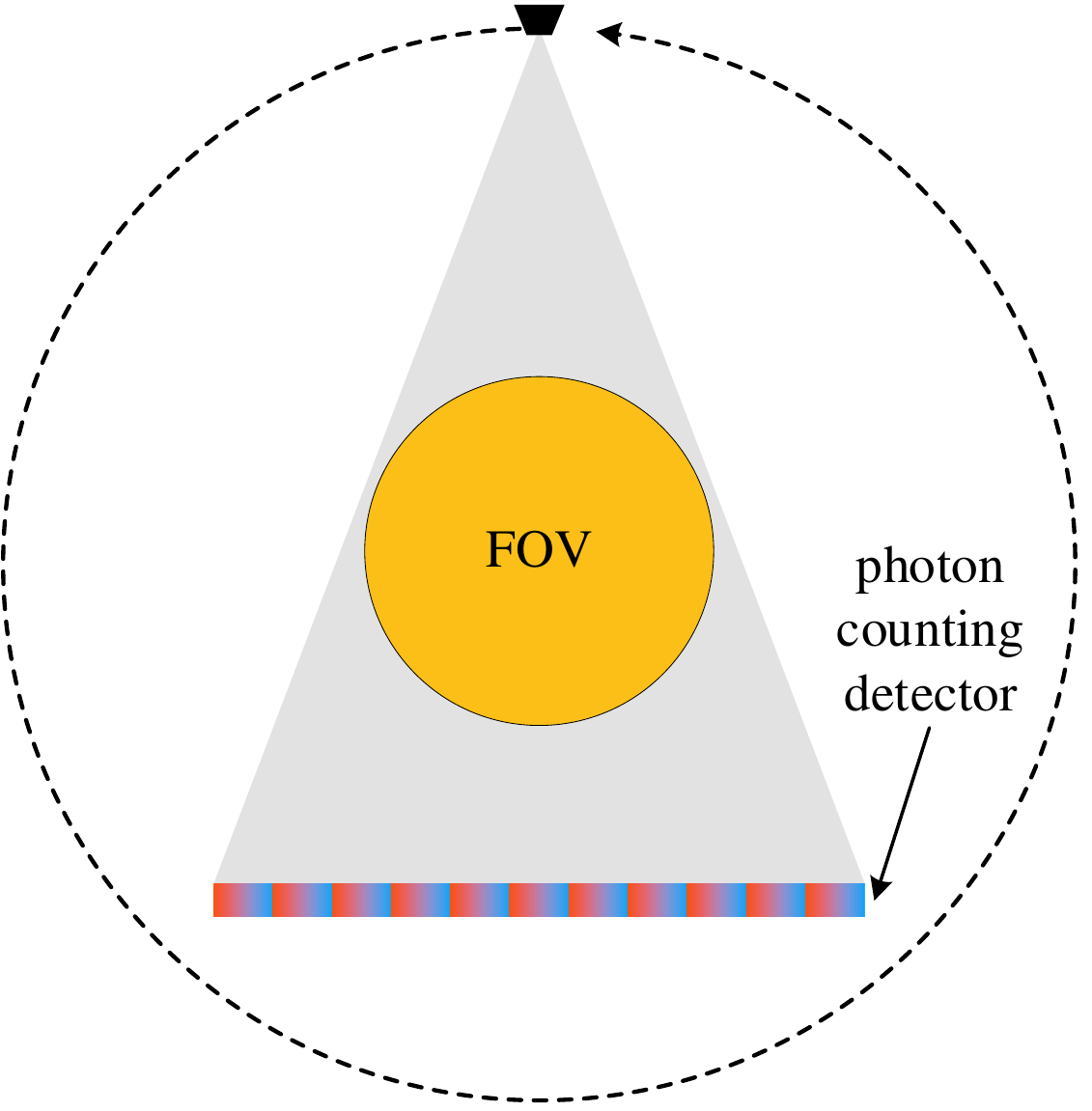}
		\label{fig1_e}}
	\subfigure[]{
		\includegraphics[width=4.5cm]{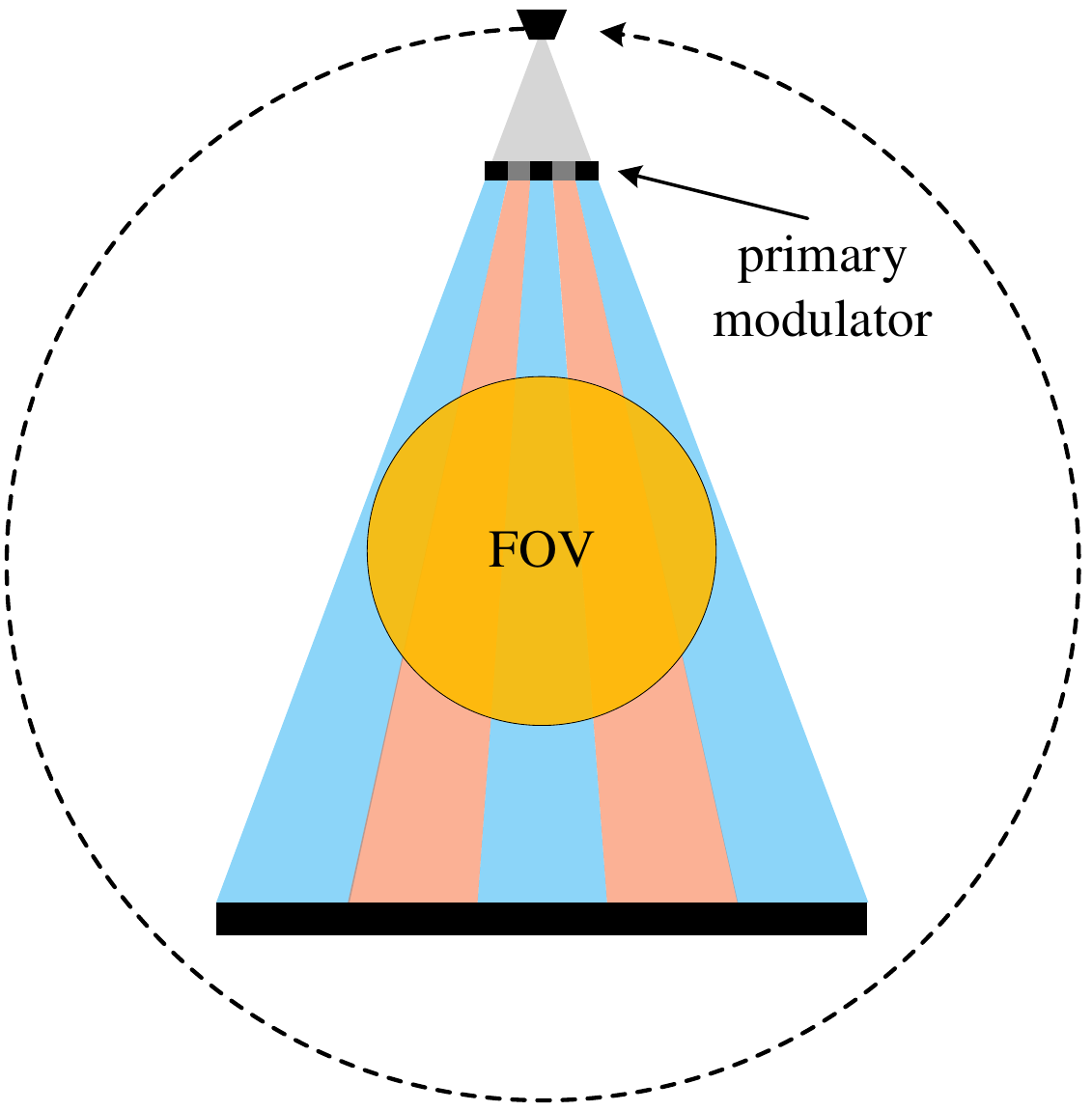}
		\label{fig1_f}} \hspace{-2cm}
	\vspace{-2.5mm}
	\caption{The schematic drawing of common scan configurations. (a) Multiple full scans configuration. (b) Dual-detector configuration. (c) Fast kVp switching configuration. (d) Dual-source configuration. (e) Photon counting detector configuration. (f) Primary modulation configuration. }
	\label{fig1}	
\end{figure}

Various scan configurations have been developed to get MSCT polychromatic projections \cite{ref23,ref24,ref25}. Figure \ref{fig1} shows schematic drawing of some common ways, including multiple full scans configuration (shown in figure \ref{fig1}(a))) \cite{ref26}, dual-detector configuration (shown in figure \ref{fig1}(b)) \cite{ref27,ref28}, fast kVp switching configuration (shown in figure \ref{fig1}(c)) \cite{ref29,ref30}, dual-source configuration (shown in figure \ref{fig1}(d)) \cite{ref31}, photon-counting detector configuration (shown in figure \ref{fig1}(e)) \cite{ref32,ref33,ref34} and primary modulation configuration (shown in figure \ref{fig1}(f)) \cite{ref35,ref36}. Among the mentioned scan configurations, the data obtained by dual-detector configuration and photon-counting detector configuration are geometrically consistent and the data obtained by other methods are geometrically inconsistent. Geometrically inconsistent means that, on the one hand, using each projection sets to reconstruct can show the same object, but, on the other hand, the paths of X-rays taken between different spectra are different because of the geometric inconsistent parameters \cite{ref37}.

Researchers usually model the process of obtaining MSCT data as nonlinear equations \cite{ref26,ref38}. Omitting scattered photons and taking MSCT basis material decomposition for example, the discrete nonlinear model of obtaining polychromatic projections is
\begin{equation}
\label{eq-polyproj}
p_{k,L} = -\ln\sum_{\omega=1}^{\Omega_k}s_{k,\omega}\delta e^{-\sum_{m=1}^{M}\theta_{m,\omega}q_{m,L}},\quad q_{m,L} = {\bf R}_L {\bf f}_m,
\end{equation}
where $p_{k,L}$ denotes the polychromatic projection of the $k$-th spectrum under the X-ray path $L$, $L\in \zeta_k$. $\zeta_k$ represents the X-ray path set of the $k$-th spectrum, $k = 1,2,\cdots,K$ and $K$ is the total number of spectra. The valid energy range of $k$-th normalized effective spectrum is equally divided into $\Omega_k$ intervals and the length of each interval is $\delta$. $s_{k,\omega}$ describes the sampling value of the $k$-th normalized effective spectrum at $\omega$ keV and $\sum_{\omega=1}^{\Omega_k}s_{k,\omega} = 1$. $\theta_{m,\omega}$ represents the sampling value of the mass attenuation coefficient (MAC) of the $m$-th basis material in $\omega$ keV interval, $m = 1,2,\cdots,M$ and $M$ is the total number of basis materials. $q_{m,L}$, {\it i.e.}, the so-called basis material projection, is line integral of the $m$-th basis material along the X-ray path $L$. ${\bf f}_m = (f_{m,1},f_{m,2},\cdots,f_{m,J})^\top$ denotes the discrete form of the $m$-th basis material density image, where $f_{m,j}$ is the sampling value of ${\bf f}_m$ at the $j$-th pixel. $(\bullet)^\top$ represents transpose and $J$ is the total number of image pixels. ${\bf R}_L=(r_{L,1},r_{L,2},\cdots,r_{L,J})$ is called the projection operator corresponding to the X-ray path $L$, where $r_{L,j}$ represents the contribution of the $j$-th pixel to the X-ray path $L$. In this paper, $s_{k,\omega}$ and $\theta_{m,\omega}$ are assumed to be known. The estimation of $s_{k,\omega}$ and the measurement of $\theta_{m,\omega}$ can be referred to \cite{ref39,ref40,ref41,ref42}.

\begin{figure}[htbp]
	\centering
	\includegraphics[width=2.7cm]{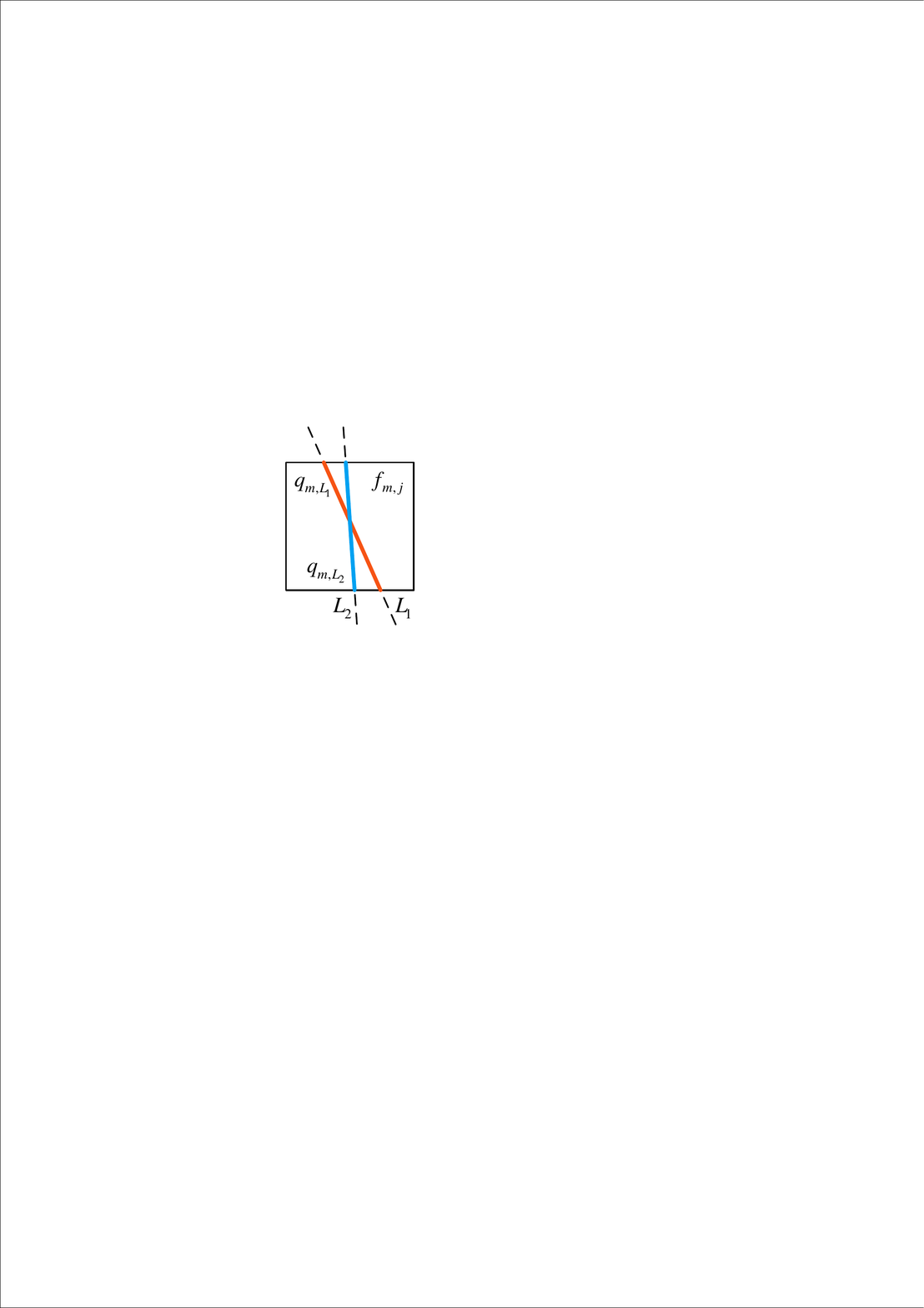} \hspace{5mm}
	\vspace{-2.5mm}
	\caption{ Geometric illustration of the geometrically inconsistent rays. The two dotted lines represent $L_1$ and $L_2$ respectively and the intersection lengths with the pixel $f_{m,j}$ are $q_{m,L_1}$ (the red line) and $q_{m,L_2}$ (the blue line) respectively.}
	\label{fig2}
\end{figure}

Dual spectral CT (DSCT) data two basis material decomposition is a typical case of MSCT basis material decomposition. The process to get DSCT polychromatic projection is
\begin{equation}
	\left\{
	\begin{aligned}
		p_{1,L_1} &= -\ln\sum_{\omega=1}^{\Omega_1}s_{1,\omega}\delta e^{-\theta_{1,\omega}q_{1,L_1}-\theta_{2,\omega}q_{2,L_1}}, \\
		p_{2,L_2} &= -\ln\sum_{\omega=1}^{\Omega_1}s_{2,\omega}\delta e^{-\theta_{1,\omega}q_{1,L_2}-\theta_{2,\omega}q_{2,L_2}}.
	\end{aligned}
	\right.
	\label{eq-DSCT}
\end{equation}
In the case of geometric consistency, $L_1$ and $L_2$ coincide. At this time, the nonlinear equations (\ref{eq-DSCT}) contains two unknowns and is a well-posed problem. However, the great majority of measured data are geometrically inconsistent. Figure \ref{fig2} shows the so-called geometrically inconsistent rays. It is clear that $L_1$ and $L_2$ pass through the same pixel $f_{m,j}$, but $q_{m,L_1} \neq q_{m,L_2} (m=1,2)$. In this case, the nonlinear equations (\ref{eq-DSCT}) contains four unknowns and is an underdetermined problem, which called as mismatched nonlinear equations in this paper.

MSCT basis material decomposition can be summarized as reconstructing basis material density image ${\bf f}_m$ from measured data $p_{k}$ by inversion of the nonlinear equations (\ref{eq-polyproj}). Works exist on investigating mapping methods, deep learning methods or iterative methods for the inversion of the nonlinear equations. The mapping methods \cite{ref43,ref44,ref45,ref46,ref47,ref48,ref49} establish the mapping relation in advance between the polychromatic projections and the basis material images. The decomposition accuracy is limited by the mapping relationship, which means that high accuracy leads to high solution complexity and high noise sensitivity. CT reconstruction based on the deep learning is a hot issue \cite{ref50}. Deep learning methods \cite{ref51,ref52,ref53,ref54} have been used in MSCT reconstruction and basis material decomposition. However, in many cases, CT data for training are difficult to obtain, especially industrial CT data. The iterative methods are most commonly used to solve the nonlinear equations. Researchers use statistical model or algebraic model to construct different iterative schemes and obtain high precision solutions by gradual correction. The iterative methods based on the statistical model take the noise distribution into account and can obtain high signal-to-noise ratio results in case of high noise \cite{ref55,ref56,ref57}. The iterative methods based on the algebraic model either invert the nonlinear model directly \cite{ref26,ref58}, or convert the nonlinear model into linear model and then solve it by linear methods \cite{ref59,ref60,ref61,ref62}. On the foundation of statistical model or the algebraic model, some researchers introduce prior information and propose optimization problems to further improve accuracy of solution \cite{ref63,ref64,ref65,ref66,ref67,ref68}.

Only focusing on the solution of the nonlinear equations, this paper summarizes most iterative methods into three steps:
\begin{description}
	\item[Step1] \textbf{Decomposition} In this step, the nonlinear equations are solved to get the basis material projection $q_{m}^{(n+1)}$.
	\item[Step2] \textbf{Reconstruction} In this step, traditional reconstruction methods, such as ART, FBP, {\it etc.}, are performed to reconstruct the basis material image ${\bf f}_{m}^{(n+1)}$ from $q_{m}^{(n+1)}$.
	\item[Step3] \textbf{Update} In this step, the new ${\bf f}_{m}^{(n+1)}$ is used to update the new polychromatic projection $p_k^{(n+1)}$ and get the new nonlinear equations. 
\end{description}

For most iterative methods, the latter two steps, {\it i.e.}, the reconstruction step and the update step, are the same. Difference appears in the decomposition step, which is the inversion of the nonlinear equations. Alvarez uses Newton-Raphson method to solve nonlinear equations \cite{ref9}. The Alvarez's method can calculate accurate solutions and have fast convergence speed for noise-free and geometrically consistent data. However, it has poor noise resistance and can not apply to geometrically inconsistent data. Several iterative methods are introduced below, which can deal with data collected with almost all configurations. Our team extended the classic ART method (E-ART) to solve the nonlinear system \cite{ref59}. The E-ART method fits almost all scanning configurations because it iterates ray by ray. High precision basis material images are reconstructed with the E-ART method. Chen proposed the ASD-NC-POCS method in 2017 \cite{ref26}. Considering the solution of the fidelity term (or the data term), the ASD-NC-POCS method combines the spectrum and the attenuation coefficient to obtain the linear part of the nonlinear model, and uses POCS to solve it. In 2021, Chen modified intercept of the ASD-NC-POCS method and developed a non-convex primal-dual (NCPD) method to solve a non-convex optimization program based on the nonlinear model \cite{ref58}. The NCPD method can yield accurate results and can be applied to deal with data collected with non-standard configurations. The above three methods can solve matched or mismatched nonlinear equations and obtain high-quality basis material images, but they have slow convergence speed.

How to solve the nonlinear equations accurately and quickly is still a hot issue. This paper proposes a general iterative method to invert the nonlinear equations and develops Schmidt orthogonalization to improve convergence speed. The method is hereafter referred to as Schmidt orthogonal modification algorithm (SOMA). For the convenience of expression, the mark of X-ray path $L$ and the length of interval $\delta$ in (\ref{eq-polyproj}) are omitted in the rest of this article.

The remainder of the article is organized as follows. In section \ref{sect2}, the principle and detailed implementation of the proposed method will be shown. In section \ref{sect3}, the simulation MSCT data and real MSCT data experiments will verify some characteristics of the proposed method. The discussion will be given in section \ref{sect4}.

\section{Method}
\label{sect2}
This section first introduces the main idea of the proposed method in the case of matched nonlinear equations, then gives the idea in the case of mismatched nonlinear equations and the general iteration scheme. Finally, the pseudo-code and the detailed implementation of applying the proposed method to MSCT basis material decomposition are shown.

\subsection{Main idea}
\label{sect2-1}
Assuming that the nonlinear equations is 
\begin{equation}
	\label{eq-gx}
	\left\{
	\begin{aligned}
		&G_1({\bf x}) = p_1 \\
		&\qquad \vdots \\
		&G_K({\bf x}) = p_K,
	\end{aligned}
	\right.
\end{equation}
where $G_k$ represents the $k$-th nonlinear equation and ${\bf x} = (x_1,x_2,\cdots,x_M)^\top$ is the unknowns. Let ${\bf x}_0$ represents the initial iteration point. Performing the first-order Taylor expansion at ${\bf x}_0$, a linear equations ${\bf A}{\bf x} = {\bf b}$ can be obtained, where
\begin{equation}
	\label{eq-Ak}
	{\bf A}_k = \nabla G_k({\bf x}_0) = (\frac{\partial G_k}{\partial x_1}({\bf x}_0),\cdots,\frac{\partial G_k}{\partial x_M}({\bf x}_0)),
\end{equation}
\begin{equation}
	\label{eq-bk}
	b_k = p_k + {\bf A}_k{\bf x}_0 - G_k({\bf x}_0).
\end{equation}
Obviously, the normal direction of the $k$-th tangent plane $H_k$ is ${\bf g}_k = {\bf A}_k^\top$.

\begin{figure}[htbp]
	\centering
	\includegraphics[width=8.5cm]{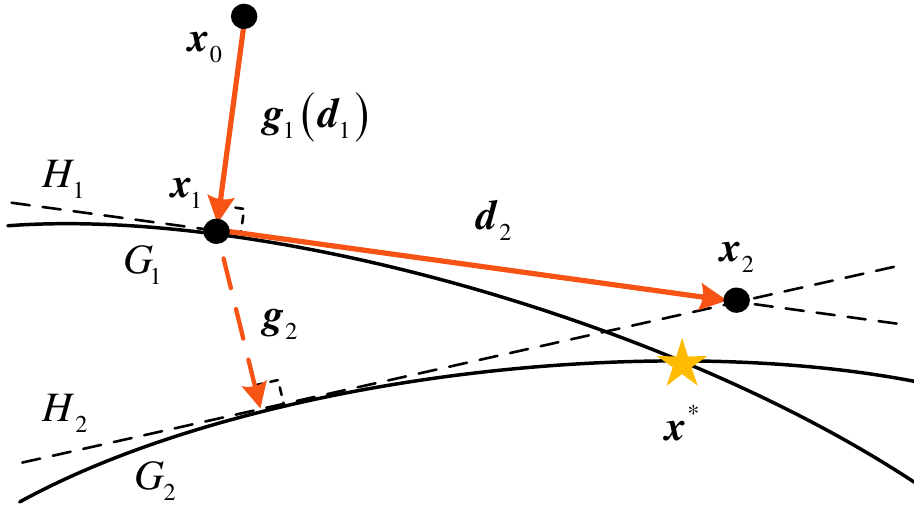} \hspace{-1cm}
	\vspace{-2.5mm}
	\caption{ The simplified geometric illustration of the proposed method. The true solution ${\bf x}^*$ is represented by the yellow pentagram. $G_1$ and $G_2$ are two surfaces, and $H_1$ and $H_2$ are the corresponding tangent planes obtained by the first-order Taylor expansion.}
	\label{fig3}
\end{figure}

A simplified geometric illustration is shown in figure \ref{fig3}. The most common method, {\it i.e.}, the gradient descent method, searches the next iteration point along the normal directions ${\bf g}_k$. However, there are usually repeated information between ${\bf g}_1$ and ${\bf g}_2$, which is easily observed in figure \ref{fig3}. Thus, the search path will be zigzag, which leads to slow convergence speed \cite{ref62}. 
 
The proposed method let the search directions as orthogonal as possible, which can reduce "redundancy" and "repeat" between the normal directions by the Schmidt orthogonalization. The orthogonal direction ${\bf d}_k$ can be obtained by the Schmidt orthogonalization
\begin{equation}
	\label{eq-dk}
	{\bf d}_k = {\bf P}_k{\bf g}_k,
\end{equation}
where ${\bf P}_k$ called orthogonal correction matrix and the initial ${\bf P}_1$ is set as the unit matrix. The update formula of ${\bf P}_k$ can be obtained by recursion
\begin{equation}
	\label{eq-pk}
	{\bf P}_k = {\bf P}_{k-1}-\frac{{\bf d}_{k-1}{\bf d}_{k-1}^\top}{{\bf d}_{k-1}^\top{\bf d}_{k-1}+\epsilon},
\end{equation}
where $\epsilon$ is generally set to a very small value in order to ensure that the denominator is not zero. As the detailed derivation of the recursion process can be found in \fullref{Appendix A}.

The next iteration point ${\bf x}_k$ can be searched for along ${\bf d}_k$. ${\bf x}_k$ is the minimum of a linear manifold $\{{\bf x}|{\bf x} = {\bf x}_{k-1} + \alpha_k{\bf d}_k\}$ spanned at ${\bf x}_{k-1}$ by $H_k$ and ${\bf d}_1,\cdots,{\bf d}_k$, where $\alpha_k$ is the optimal step size. The problem on the linear manifold is transformed into finding $\alpha_k$ satisfied
\begin{equation}
	\label{eq-min-alpha}
	\alpha_k = \min_\alpha \Vert{\bf A}_k({\bf x}_{k-1} + \alpha{\bf d}_k) - b_k\Vert^2_2.
\end{equation}
Let the partial derivative with respect to $\alpha$ be equal to 0 and it is easy to calculate solution
\begin{equation}
	\label{eq-alpha}
	\alpha_k=\frac{{\bf d}_k^\top{\bf g}_k ( b_k - {\bf A}_k{\bf x}_{k-1}) }{ {\bf d}_k^\top {\bf g}_k {\bf g}_k^\top {\bf d}_k} = \frac{b_k - {\bf A}_k{\bf x}_{k-1} }{{\bf g}_k^\top {\bf d}_k}.
\end{equation}

A simple convergence proof of the SOMA method is given in \fullref{Appendix B}.

\subsection{The SOMA method for the mismatched nonlinear equations}
\label{sect2-2}
For the mismatched nonlinear equations, only one equation can be obtained accurately and the constant term of other equations are unknown. Assume the first nonlinear equation is accurate and the corresponding linear equation is ${\bf A}_1{\bf x} = b_1$. Referring to \cite{ref38}, other linear equations can be shown as follow
\begin{equation}
	\label{eq-mis-linear}
	{\bf A}_k{\bf x} = \tilde{b}_k,
\end{equation}
where ${\bf A}_k$ is defined by formula (\ref{eq-Ak}) and $\tilde{b}_k = $ represents the unknown intercept
\begin{equation}
	\tilde{b}_k = \tilde{p}_k + {\bf A}_k{\bf x}_0 - G_k({\bf x}_0),
\end{equation}
in which $\tilde{p}_k$ represents the unknown constant term.

\begin{figure}[htbp]
	\centering
	\includegraphics[width=8.5cm]{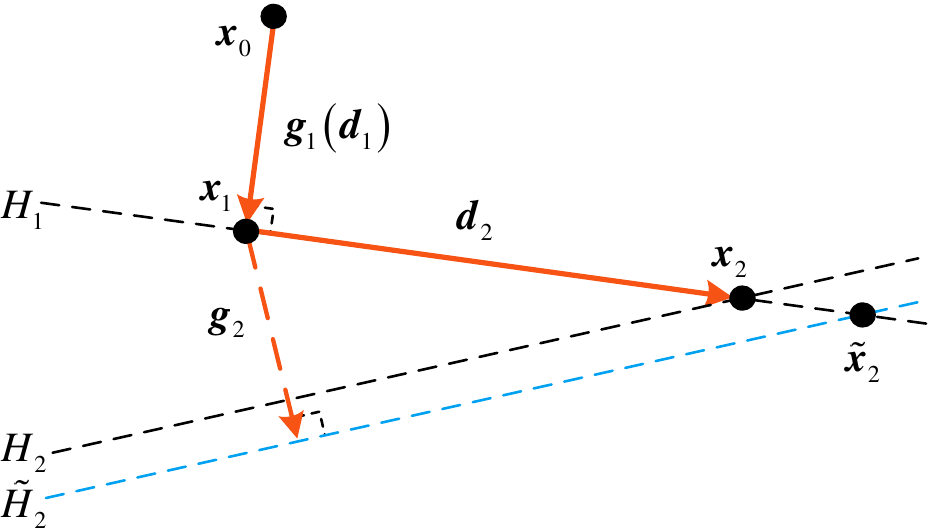} \hspace{-1cm}
	\vspace{-2.5mm}
	\caption{ The simplified geometric illustration in the case of the mismatched nonlinear equations. The black dotted lines $H_k$ are the accurate tangent planes. The blue dotted line $\tilde{H}_2$ represents the tangent plane with unknown intercept. The curves $G_k$ are omitted for showing clearly.}
	\label{fig4}
\end{figure}

Figure \ref{fig4} shows the simplified geometric illustration of the mismatched nonlinear equations. The black dotted lines are the accurate tangent planes $H_k$. However, because of the unknown intercept $\tilde{b}_2$, $H_2$ could not be get. Only the tangent plane $\tilde{H}_2$ with the unknown intercept can be obtained, which is represented by the blue dotted line. Obviously, if still using the step size $\alpha_k$ calculated by the formula (\ref{eq-alpha}), the wrong iteration point $\tilde{{\bf x}}_{2}$ will be obtained.

One way to get the true optimal step size is the one-dimensional search method. Or, there is another simple way that adding an step size relaxation factor $\beta$
\begin{equation}
	\label{eq-xnk-mis}
	{\bf x}_k = {\bf x}_{k-1} + \beta\alpha_k{\bf d}_k.
\end{equation}

Therefore, the general iteration scheme of the proposed method as follows:
\begin{algorithm}[htbp]
	\setstretch{1.15}
	\caption{The general iteration scheme of the SOMA method.}
	\label{alg1}
	\LinesNumbered
	\textbf{initialize:} assign ${\bf x}_0$, $\beta$ and $\epsilon$ with some initial values\\
	\While{not satisfying the stopping criterion}{
		Perform the first-order Taylor expansion at ${\bf x}_0$ and get the linear equations \\
		${\bf P}_1 = {\bf I}$ \\
		${\bf g}_k = {\bf A}_k^\top$ \\
		\For{$k=1$ to $K$}{
		${\bf d}_k = {\bf P}_{k}{\bf g}_k $ \\
		$\alpha_k=\frac{ b_k -{\bf A}_k {\bf x}_{k-1} }{{\bf g}_k^\top {\bf d}_k}$ \\
		${\bf x}_k = {\bf x}_{k-1} + \beta\alpha_k{\bf d}_k$ \\
		${\bf P}_{k+1} = {\bf P}_{k}-\frac{{\bf d}_k{\bf d}_k^\top}{{\bf d}_k^\top{\bf d}_k+\epsilon}$
	}
	${\bf x}_0 = {\bf x}_K$
	}
\end{algorithm}

A variant of the proposed method is introduced in \fullref{Appendix C}. Although its convergence proof has not been given, its effect is better than the method in this section according to the current research.

\subsection{Implementation in MSCT basis material decomposition}
\label{sect2-3}
In section \ref{sect1}, this paper summarizes most iterative methods into three steps. For MSCT basis material decomposition, the proposed method is applied to the decomposition step, {\it i.e.}
\begin{description}
	\item[Step1] \textbf{Decomposition} Using the SOMA method to solve the nonlinear equations and get the basis material projection $q_{m}^{(n+1)}$.
	\item[Step2] \textbf{Reconstruction} Applying the traditional reconstruction methods, such as ART, FBP, {\it etc.}, to reconstruct the basis material image ${\bf f}_{m}^{(n+1)}$ from $q_{m}^{(n+1)}$.
	\item[Step3] \textbf{Update} Updating the new polychromatic projection $p_k^{(n+1)}$ and get the new nonlinear equations.
\end{description}

The pseudo-code of applying the proposed method to MSCT basis material decomposition is shown in algorithm \ref{alg2}. To avoid adjusting parameters, the pseudo-code contains the adaptive step size strategy. After the pseudo-code, the detailed implementation is explained.

\begin{algorithm}[htbp]
	\setstretch{1.1}
	\caption{Pseudo-code of applying the SOMA method to MSCT basis material decomposition.}
	\label{alg2}
	\LinesNumbered
	\textbf{initialize:} set ${\bf f}_m^{(0)} = 0$, $\lambda = 0.9$, $\epsilon = 10^{-8}$, $T = 1.5$, $\beta = 0.9$ and $\beta_{red} = 0.9$\\
	\While{not satisfying the stopping criterion}{
		\If{the polychromatic projections are inconsistent}{
			estimate the unknown projection $\tilde{p}_k$ under the current path and let $p_k=\tilde{p}_k$ \\
		}
		${\bf q}^{(n)} = [{\bf R}{\bf f}_1^{(n)},{\bf R}{\bf f}_2^{(n)},\cdots,{\bf R}{\bf f}_M^{(n)}]^\top$ \\
		Perform the first-order Taylor expansion at ${\bf q}^{(n)}$ and get the linear equations ${\bf A}_k^{(n)}{\bf q} = b_k$ \\ 
		${\bf g}_k^{(n)} = {\bf A}_k^{(n)\top} =  [\frac{\Theta_{k,1}^{(n)}}{\Phi_{k}^{(n)}}, \frac{\Theta_{k,2}^{(n)}}{\Phi_{k}^{(n)}}, \cdots, \frac{\Theta_{k,M}^{(n)}}{\Phi_{k}^{(n)}}]^\top$ \\
		${\bf q}^{(n,0)} = {\bf q}^{(n)}$ \\
		${\bf P}_0^{(n)} = {\bf I}$ \\
		\For{$k=0$ to $K-1$}{
			${\bf d}_k^{(n)} = {\bf P}_{k}^{(n)}{\bf g}_k^{(n)}$ \\
			$\alpha_k^{(n)} = \frac{ b_k - {\bf A}_k {\bf x}_{k-1}}{{\bf g}_k^\top {\bf d}_k}$ \\
			${\bf q}^{(n,k+1)} = {\bf q}^{(n,k)} + \beta\alpha_k^{(n)}{\bf d}_k^{(n)}$ \\
			${\bf P}_{k+1}^{(n)} = {\bf P}_{k}^{(n)}-\frac{{\bf d}_k^{(n)}{\bf d}_k^{(n)\top}}{{\bf d}_k^{(n)\top}{\bf d}_k^{(n)}+\epsilon}$ \\
		}
		$dp = \frac{\Vert p_k-p_k^{(n,K)} \Vert_2^2}{\Vert p_k-p_k^{(n,1)} \Vert_2^2}$ \\
		$df_m = \frac{\Vert {\bf f}_m^{(n)}-{\bf R}^{-1}(q_m^{(n,K)}) \Vert_2^2}{\Vert {\bf f}_m^{(n)}-{\bf R}^{-1}(q_m^{(n,1)}) \Vert_2^2}$ \\
		\If{$dp>1$ or $df_m \ge T$}{
			${\bf q}^{(n,K)} = {\bf q}^{(n,1)}$ \\
			$\beta = \beta\cdot\beta_{red}$ \\
		}
		${\bf q}^{(n+1)} = {\bf q}^{(n,K)}$ \\
		${\bf f}_m^{(n+1)} = {\bf f}_m^{(n)} + \lambda\cdot {\bf R}^{-1}(q_m^{(n+1)}-q_m^{(n)})$ \\
	}
\end{algorithm}

Line 1 gives the initial values of parameters and the basis material images. Lines 2-25 form the loop part of the proposed method, which can be divided into three modules, using the SOMA method to perform MSCT basis material decomposition (lines 7-16), the adaptive step size strategy (lines 17-22), and updating the estimated values of the basis material images (lines 23-24). Each module is described in detail below.

Lines 3-5 are only performed when the polychromatic projections are geometrically inconsistent. There are many methods to estimate $\tilde{p}_k$, such as the nearest neighbor interpolation method and the linear interpolation method. When the sampling is dense enough, the estimated projections $\tilde{p}_k$ is a good approximation. On the other hand, $\tilde{p}_k$ is just used to calculate $\alpha$. Even if the sampling is not dense enough, the step size relaxation factor $\beta$ can reduce some errors.

Line 6 is getting the basis material projection values of the $n$-th iterations ${\bf q}^{(n)}$. Line 7 is performing the first-order Taylor expansion of (\ref{eq-polyproj}) at ${\bf q}^{(n)}$ to obtain the linear equation, where
\begin{eqnarray}
	\label{eq-theta}
	b_k &= p_k + {\bf A}_k^{(n)}{\bf q}^{(n)} - p_k^{(n)}, \\
	p_{k}^{(n)} &= -\ln\sum_{\omega=1}^{\Omega_k} s_{k,\omega} e^{-\sum_{m=1}^{M}\theta_{m,\omega}q_{m}^{(n)}}, \\
	\Theta_{k,m}^{(n)} &= \sum_{\omega=1}^{\Omega_k}s_{k,\omega}\theta_{m,\omega} e^{-\sum_{m=1}^{M}\theta_{m,\omega}q_{m}^{(n)}}, \\
	\Phi_{k}^{(n)} &= \sum_{\omega=1}^{\Omega_k}s_{k,\omega} e^{-\sum_{m=1}^{M}\theta_{m,\omega}q_{m}^{(n)}}.
\end{eqnarray}

Line 8-16 correspond to line 4-11 of algorithm \ref{alg1}, which using the SOMA method to calculate the basis material projection values of the next iteration.

The adaptive step size strategy is shown in lines 17-22 and it adjusts the step size by changing the step size relaxation factor $\beta$. The adjustment condition is $dp$ and $df_m$, in which $dp$ is the ratio of the basis material projection residual changes between $k=1$ and $k=K$, and $df_m$ is the basis material image changes. Line 17 and line 18 calculate $dp$ and $df_m$ respectively. On the one hand, the residuals should be smaller when $k=K$ (that is $dp<1$). On the other hand, the images should not change too much (that is $df_m<T$, where $T$ is the threshold set in advance). If any above conditions are not met, it indicates that the step size may be large at this time, the true solution may be missed, and the step size should be reduced. In lines 19-22, the step size is adjusted adaptively according to the above two ratios, where $\beta_{red}$ is the reduced ratio of the step relaxation factor. Line 20 shows that the results after traversing $K$ equations are not credible, so using the results at $k = 1$ as this iteration's final result.

Line 23 is obtaining the new basis material projection values and line 24 update the basis material images by the basis material projection residuals. $\lambda$ is reconstruction relaxation factor and generally speaking, $\lambda = 1$ is enough. When higher-precision reconstruction results are required, $\lambda$ can be attenuated as the iterations increase \cite{ref60}.

A simplified DSCT example is shown in \fullref{Appendix D} and it illustrates the noticeable effect of the proposed method for convergence speed clearly. The Alvarez method \cite{ref9}, the E-ART method \cite{ref59} and the NCPD method \cite{ref58} are selected as the comparison methods because they are representative iterative methods. The Alvarez method is almost the first DSCT reconstruction technique but it only can be applied to geometrically consistent data. Its result is a benchmark in section \ref{sect3-2}. The latter two methods, not only can they inverse the mismatched nonlinear equations, but also can decompose the basis material images accurately and obtain high-quality reconstruction results.

\section{Experiment}
\label{sect3}
In this section, four MSCT basis material decomposition experiments are carried out. Firstly, a numerical convergence experiment is used to verify the numerical convergence of the proposed method. Then, the noisy data experiment shows its robustness to noise. Next, the triple material data experiment studies the feasibility of multiple basis material decomposition. Finally, a real data experiment illustrates its practical value.

\subsection{Numerical convergence verification}
\label{sect3-1}
In this section, noise-free data are used to study the numerical convergence of the proposed method. First, dual-domain convergence conditions are given. Then, the experiments are designed. The results are finally shown.

\subsubsection{Dual-domain convergence condition}
\label{sect3-1-1}

Refer to \cite{ref26}, the distance of data $D_{\mathrm{data}}^{(n)}$ and the distance of images $D_{\mathrm{image}}^{(n)}$ are given as follows
\begin{eqnarray}
\label{eq-dualdomain}
D_{\mathrm{data}}^{(n)}=\sum_{k=1}^{K} \frac{(p_k - p_{k}^{(n)})^2}{(p_k)^2}, \\
D_{\mathrm{image}}^{(n)}=\sum_{m=1}^{M} \frac{(f_m^{[\mathrm{true}]} - f_{m}^{(n)})^2}{(f_m^{[\mathrm{true}]})^2},
\end{eqnarray}
where $f_m^{[\mathrm{true}]}$ represents the true values of the $m$-th basis material images. The projection domain convergence condition is $D_{\mathrm{data}}^{(n)}\rightarrow0$ when $n\rightarrow\infty$, and the image domain convergence condition is $D_{\mathrm{image}}^{(n)}\rightarrow0$ when $n\rightarrow\infty$.

\subsubsection{Experiment design}
\label{sect3-1-2}
\begin{figure}[htbp]
	\centering
	\subfigure[]{	
		\label{fig7_a}
		\includegraphics[width=4cm]{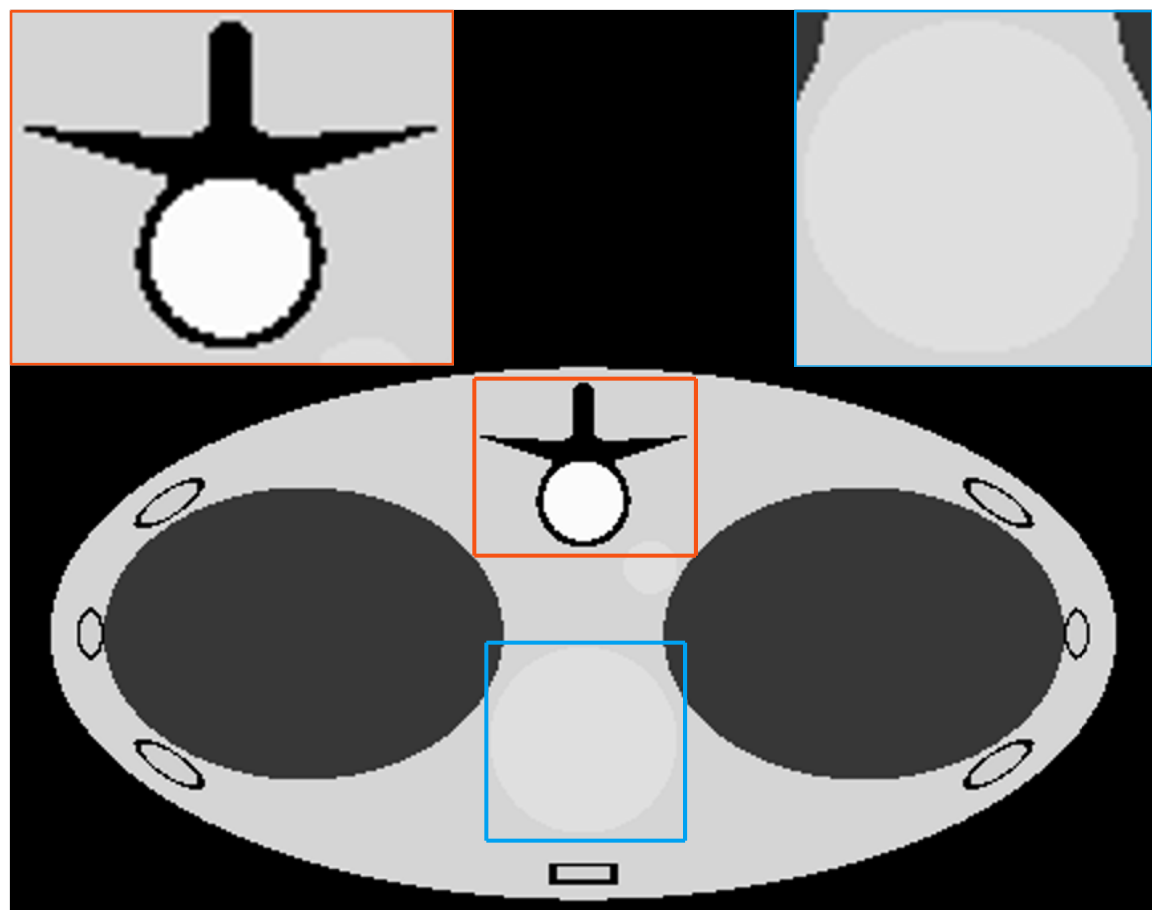}}
	\subfigure[]{	
		\label{fig7_b}	
		\includegraphics[width=4cm]{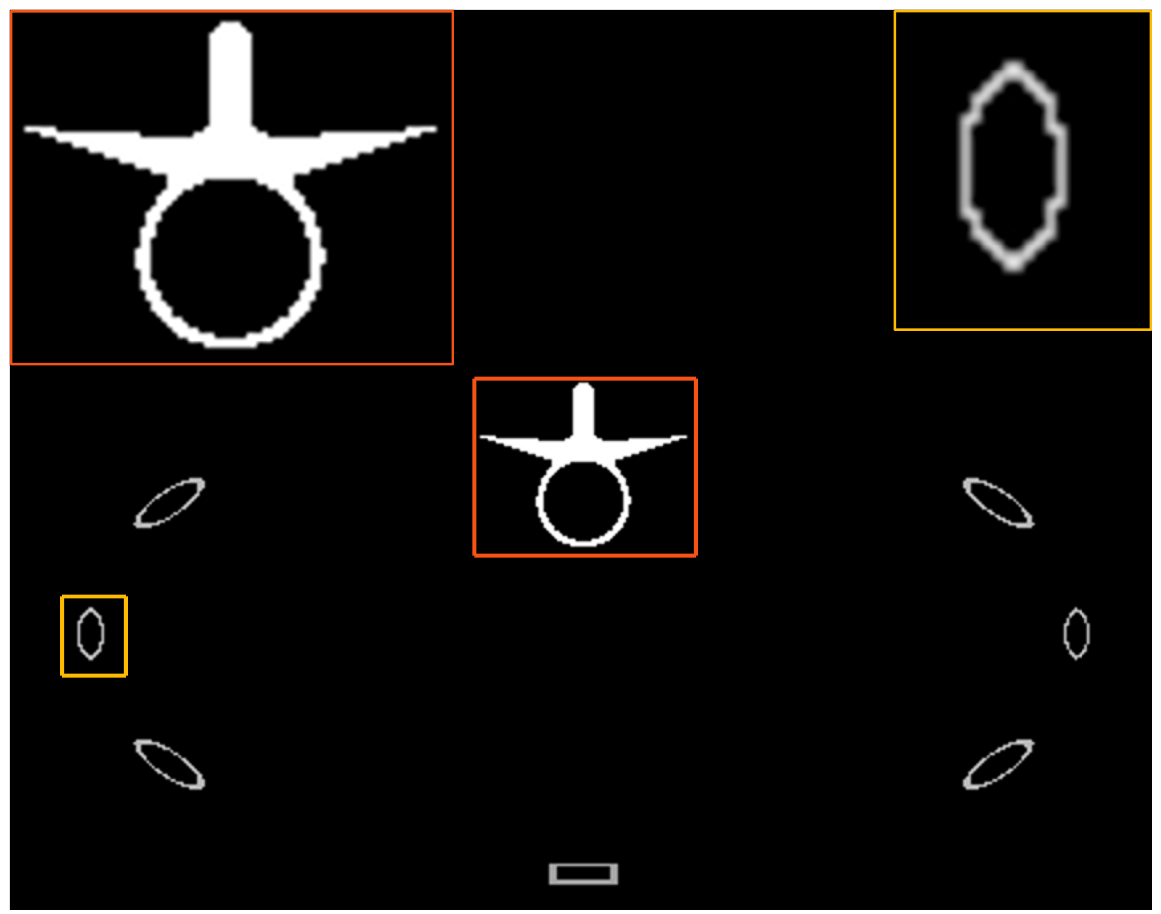}}
	\subfigure[]{	
		\label{fig7_c}
		\includegraphics[width=4cm]{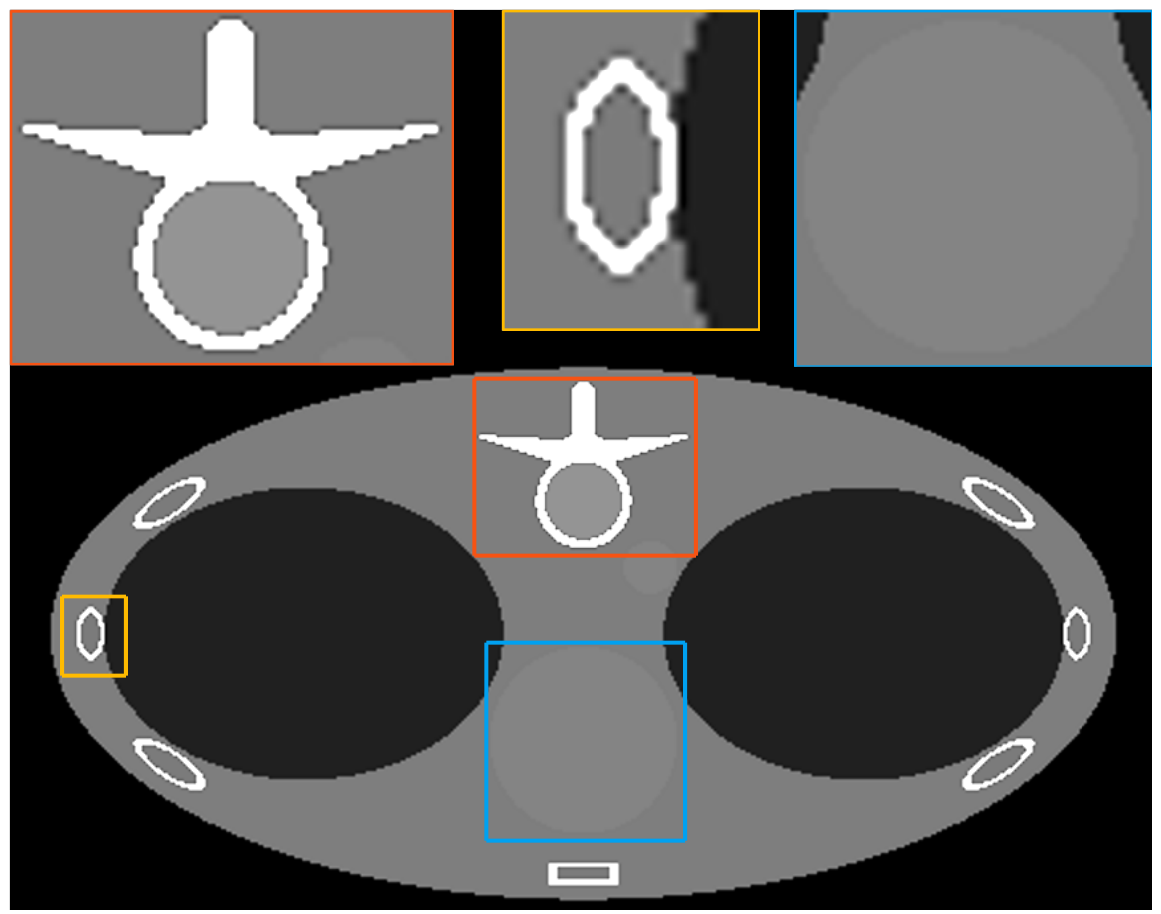}} \hspace{-2cm}
	\vspace{-2.5mm}
	
	\subfigure[]{
		\label{fig7_d}
		\includegraphics[width=7cm]{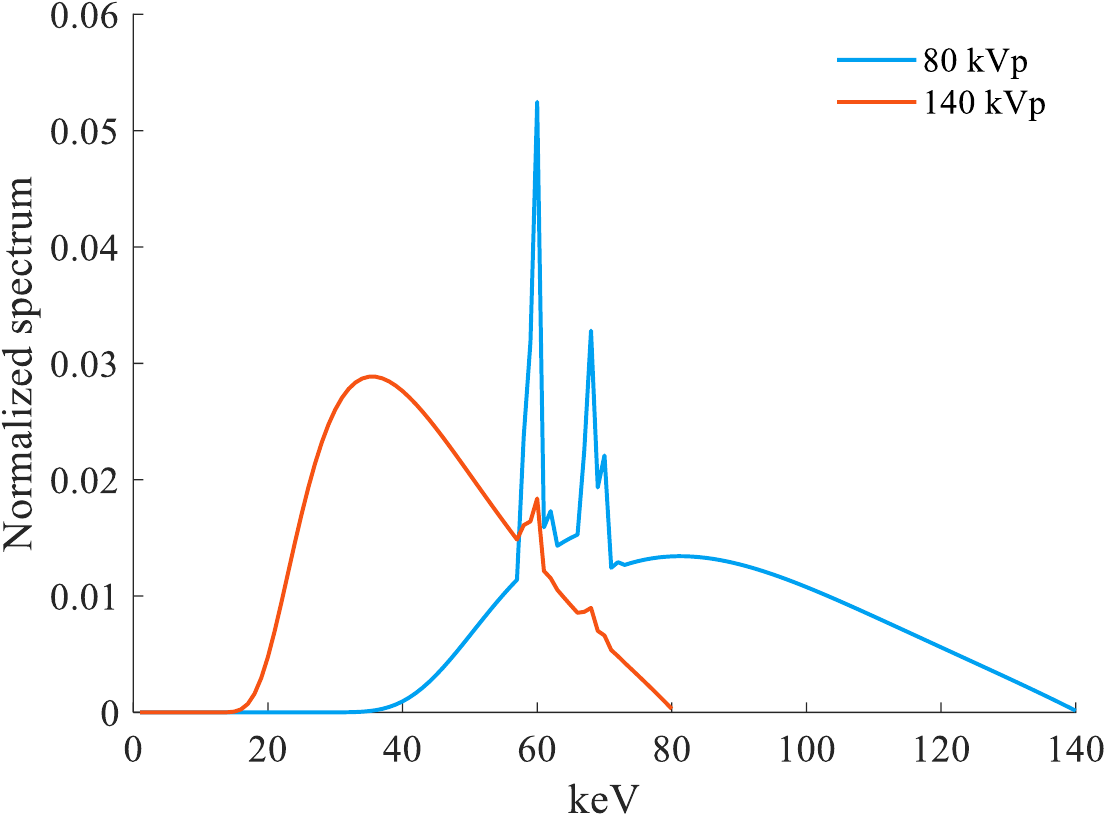}} \hspace{-2cm}
	\vspace{-2.5mm}
	\caption{Phantom and spectra used in numerical convergence experiment. (a) The water basis material image (display window: [0, 1.20]). (b) The bone basis material image (display window: [0, 1.92]). (c) The monochromatic image at 70keV (display window: [0, 0.04]). (d) Simulated spectra.}
	\label{fig7}	
\end{figure}

The phantom is a slice of the 3D FORBILD thorax phantom with resolution $512\times512$, and its geometric information and the reference density value are detailed on the website \cite{ref69}. Assuming that the phantom is composed of water and bone, the densities are 1.00 g/cm$^3$ and 1.92 g/cm$^3$ respectively. The mass attenuation coefficients of water and bone can be obtained from the National Institute of Standard Technology (NIST) website \cite{ref70}. Figure \ref{fig7}(a) and \ref{fig7}(b) shows two basis material images respectively. \ref{fig7}(c) is the monochromatic image at 70 keV and choosing this energy because of good visual effect. The spectra are generated by the open-source software Spectrum GUI \cite{ref71}, the X-ray source is GE Maxiray 125 X-ray tube, the tube voltage is set to 80 kVp and 140 kVp, and the latter is filtered by 1 mm copper, and the normalized X-ray spectra are shown in figure \ref{fig7}(d). The sampling intervals of the X-ray energy spectra and the mass attenuation coefficients are both 1 kVp.

The scanning configurations are set as follows: the distance between the X-ray source and the turntable center (SOD) is 541 mm, and the distance between the X-ray source and the detector (SDD) is 949 mm. The linear detector is composed of 960 detector cells, and the size of each cell is 1.25 mm. Using the above parameters, the diameter of the view field is approximately 579.28 mm. Simulate the process of obtaining polychromatic projections by formula (\ref{eq-polyproj}). Using multiple full scans, a total of 720 projections are collected under each spectrum for the single full scanning. When simulating the process of generating geometrically inconsistent data, the initial angle of the single full scanning differs by 0.25$^{\circ}$. The computer used in the experiment is equipped with a 2.40 GHz Intel Xeon E5-2620 six-core CPU and an NVIDIA Quadro K2200 graphics card. The size of the reconstructed estimated image is $512\times512$, and each pixel of the initial estimated image is set to 0.

\subsubsection{Experiment results}
\label{sect3-1-3}
\begin{figure}[htbp]
	\centering
	\subfigure[]{	
		\label{fig8_a}
		\includegraphics[width=7cm]{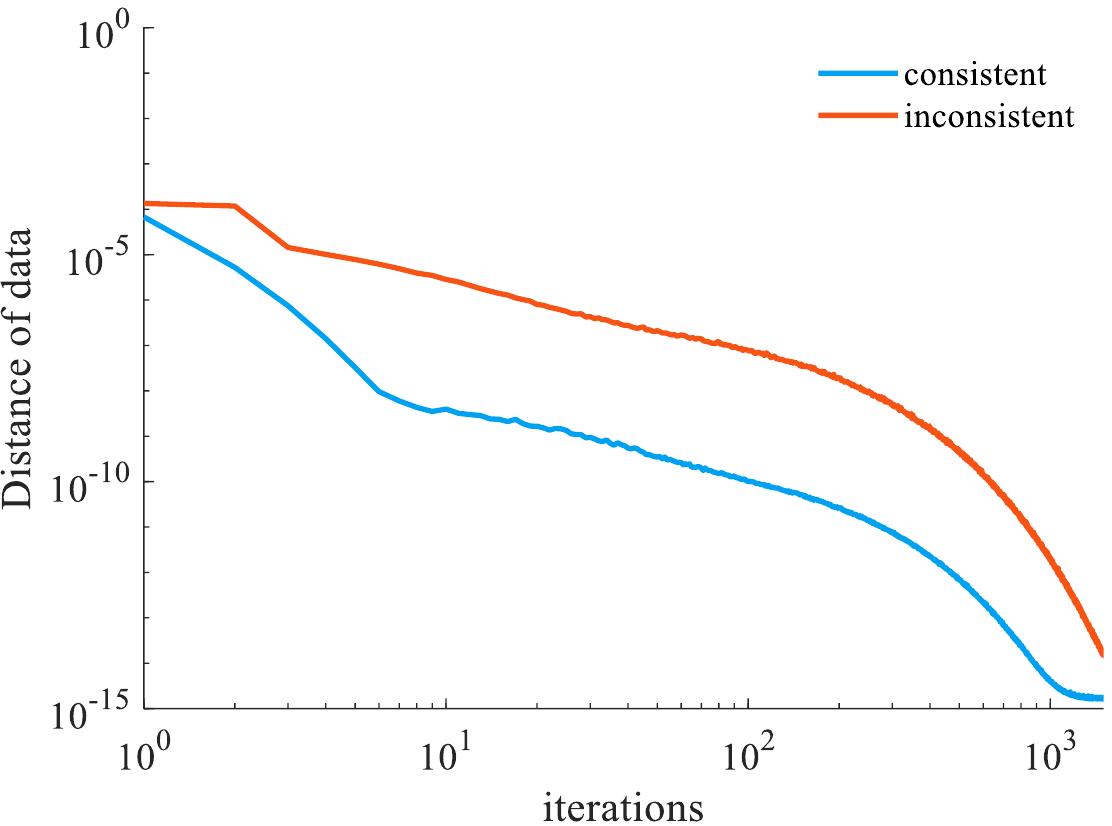}} \hspace{-5mm}
	\hspace{0.5cm}
	\subfigure[]{	
		\label{fig8_b}	
		\includegraphics[width=7cm]{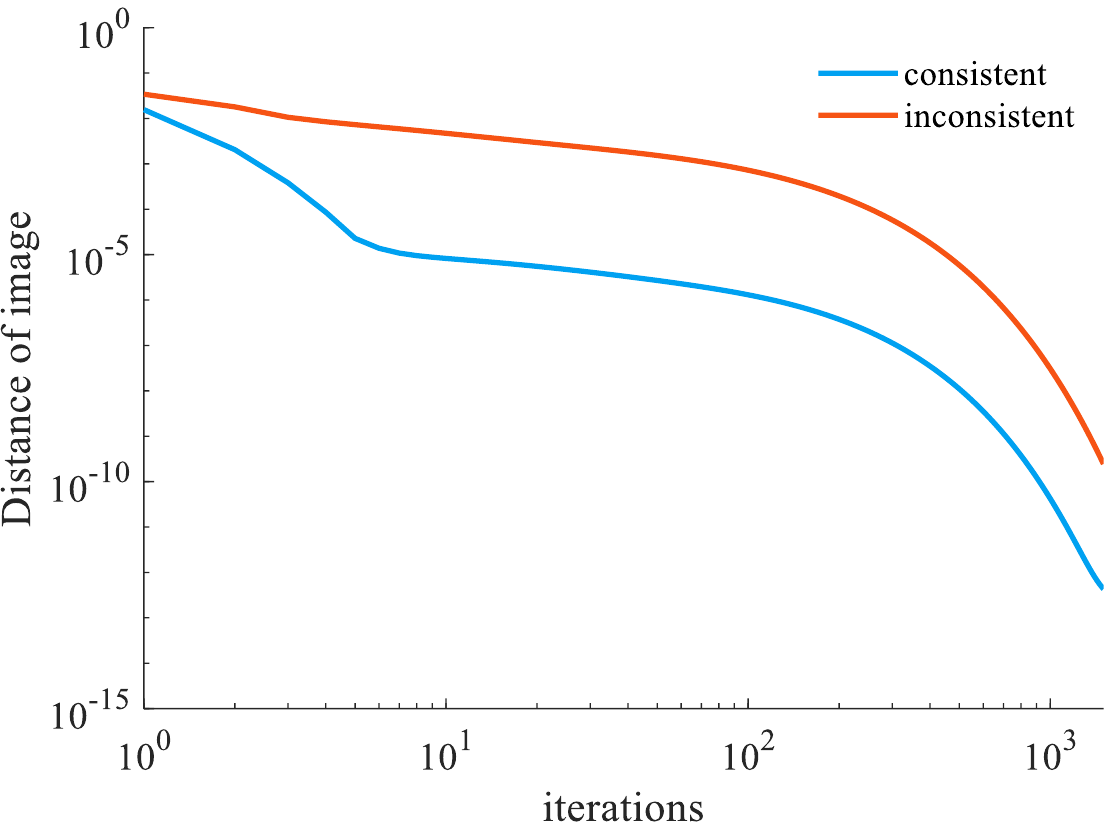}} \hspace{-2cm}
	\vspace{-2.5mm}
	\caption{Dual-domain convergence curves. (a) The convergence curves of the distance of data. (b) The convergence curves of the distance of images.}
	\label{fig8}	
\end{figure}

Using the proposed method to decompose and reconstruct the geometrically consistent and inconsistent data. For consistent data, the parameter settings are $\beta=1.0$, $\lambda=1$, $\epsilon=10^{-8}$, $\eta = 1.5$, $\kappa = 0.95$ and $\alpha = 1.0$. For inconsistent data, except $\beta = 0.5$, other parameters are the same as consistent data. Using all equation information and updating once in the image domain is considered as an iteration. We compute $D_{\mathrm{data}}^{(n)}$ and $D_{\mathrm{image}}^{(n)}$ and plot them as functions of iteration $n$ in figure \ref{fig8}. Both of them decrease with the increase of the iterations, reach the float accuracy and have a continued downward trend. Hence, the results indicate that the proposed method satisfies the dual-domain convergence conditions and is numerically convergent.

\begin{figure}[htbp]
	\centering	
	\includegraphics[width=13cm]{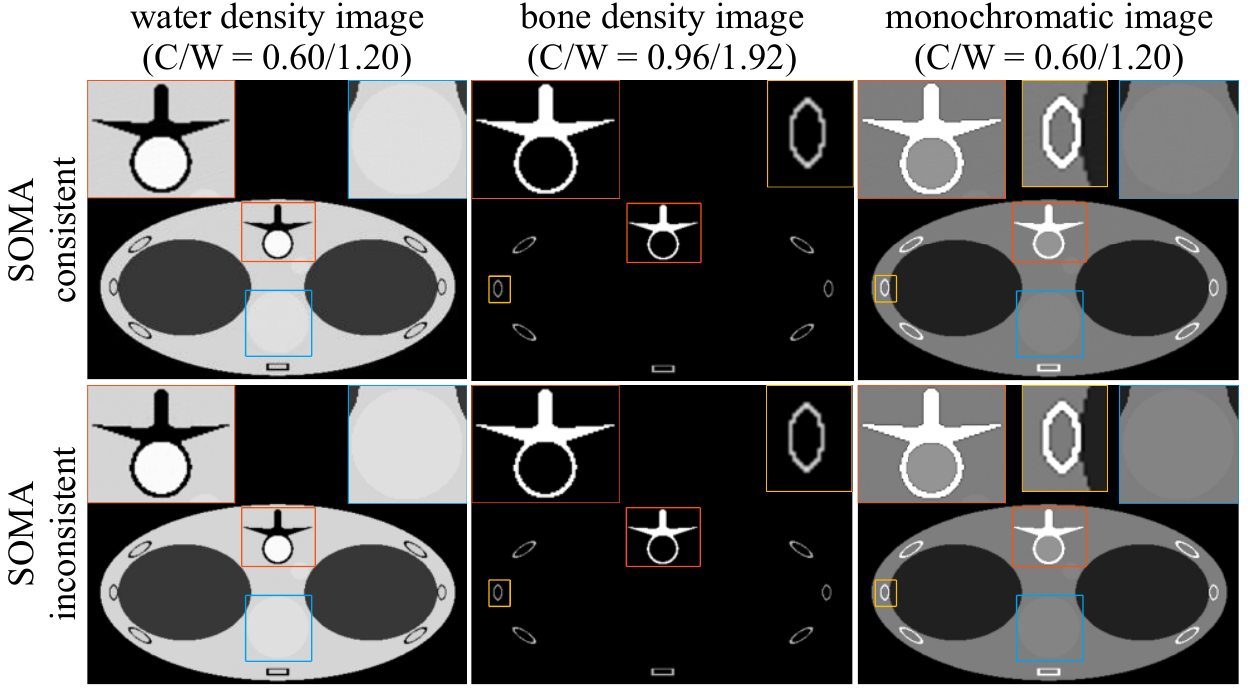} \hspace{-2.3cm}
	\vspace{-2.5mm}
	\caption{Reconstruction results of noise-free data with the proposed method when $D_{\mathrm{image}}^{(n)}<10^{-3}$. The first row shows the results from consistent data after 3 iterations. The second row shows the results from inconsistent data after 78 iterations.}
	\label{fig9}	
\end{figure}

Figure \ref{fig9} shows the reconstruction results when $D_{\mathrm{image}}^{(n)}<10^{-3}$. The results show that the proposed method can reconstruct the basis material image, and there is no visual difference between the reconstruction result and the phantom.

\subsection{Noisy data experiment}
\label{sect3-2}
The phantom, spectra, and scan configurations used in this section are the same as those in the previous section. Noisy polychromatic projections are simulated by adding Poisson noise with an initial number of $10^5$ photons to the noise-free data to study the robustness of the proposed method to noise. The Alvarez's method is applied to geometrically consistent data and the result is a benchmark. The E-ART method and the NCPD method are selected as the comparison algorithm. In order to more fairly verify the solution of the fidelity term, the parameter of the regular term of the NCPD method is set to 0. In addition, the E-ART method has only one parameter, the image domain relaxation factor $\lambda$, thus for the three methods, $\lambda$ is set to 1 in all experiments below.

\begin{figure}[]
	\centering	
	\includegraphics[width=13cm]{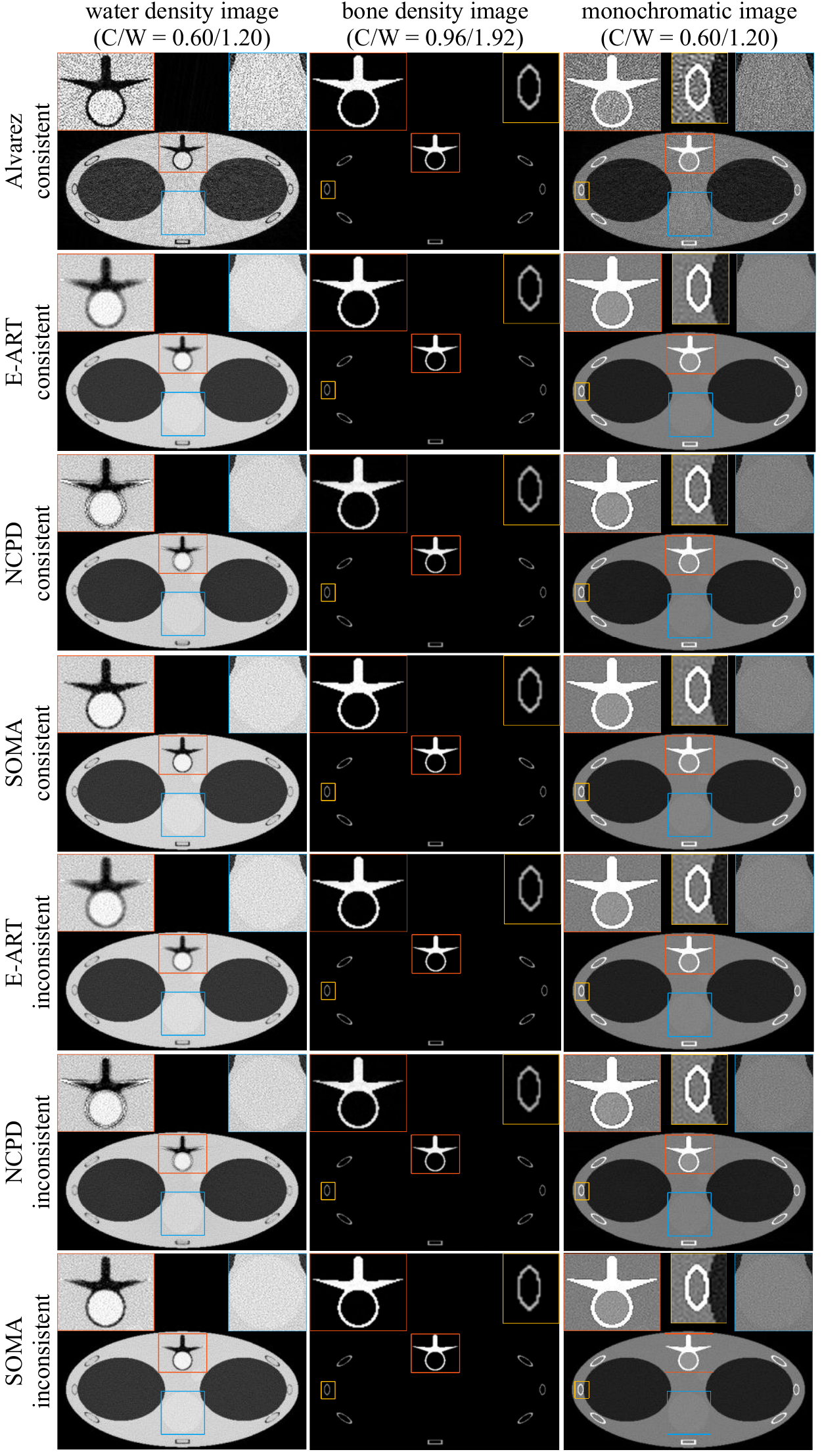} \hspace{-2.3cm}
	\vspace{-2.5 mm}
	\caption{Reconstruction results of noisy data after 30 iterations.}
	\label{fig10}	
\end{figure}

Generally speaking, it is hoped that the method needs as few iterations as possible. For this reason, we set a stopping criterion as stopping the iteration at 30 iterations in this experiment. For geometrically consistent and inconsistent data, the parameter of the proposed method are all set as $\beta=0.05$, $\epsilon=10^{-8}$, $\eta = 1.5$, $\kappa = 0.95$ and $\alpha = 0.99$. Figure \ref{fig10} shows the reconstruction results of the E-ART method, the NCPD method and the proposed method when the stopping criterion is met.

\begin{figure}[htbp]
	\centering
	\subfigure[]{	
		\label{fig11_a}
		\includegraphics[width=7.7cm]{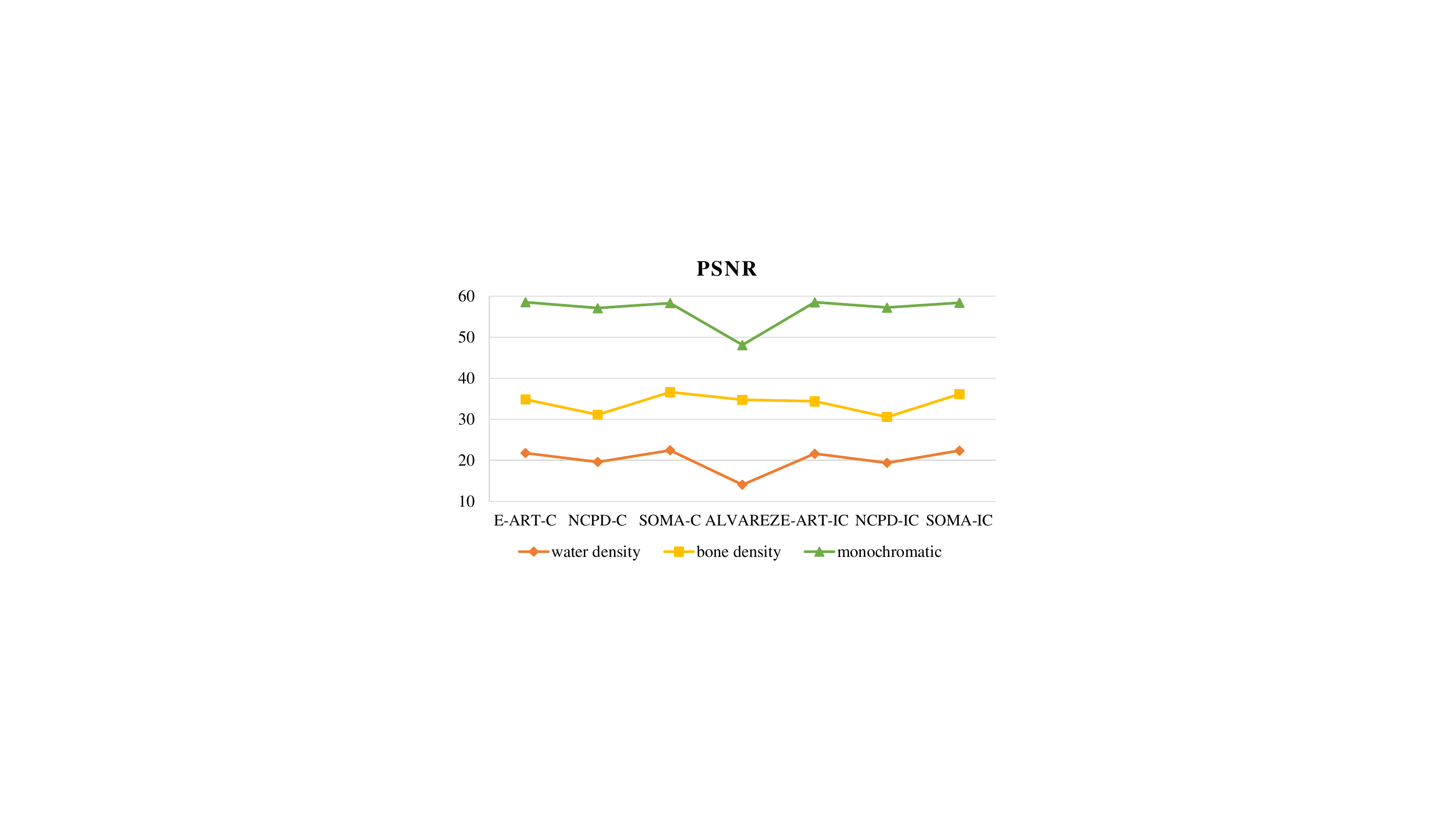}} \hspace{-2mm}
	\subfigure[]{	
		\label{fig11_b}	
		\includegraphics[width=7.7cm]{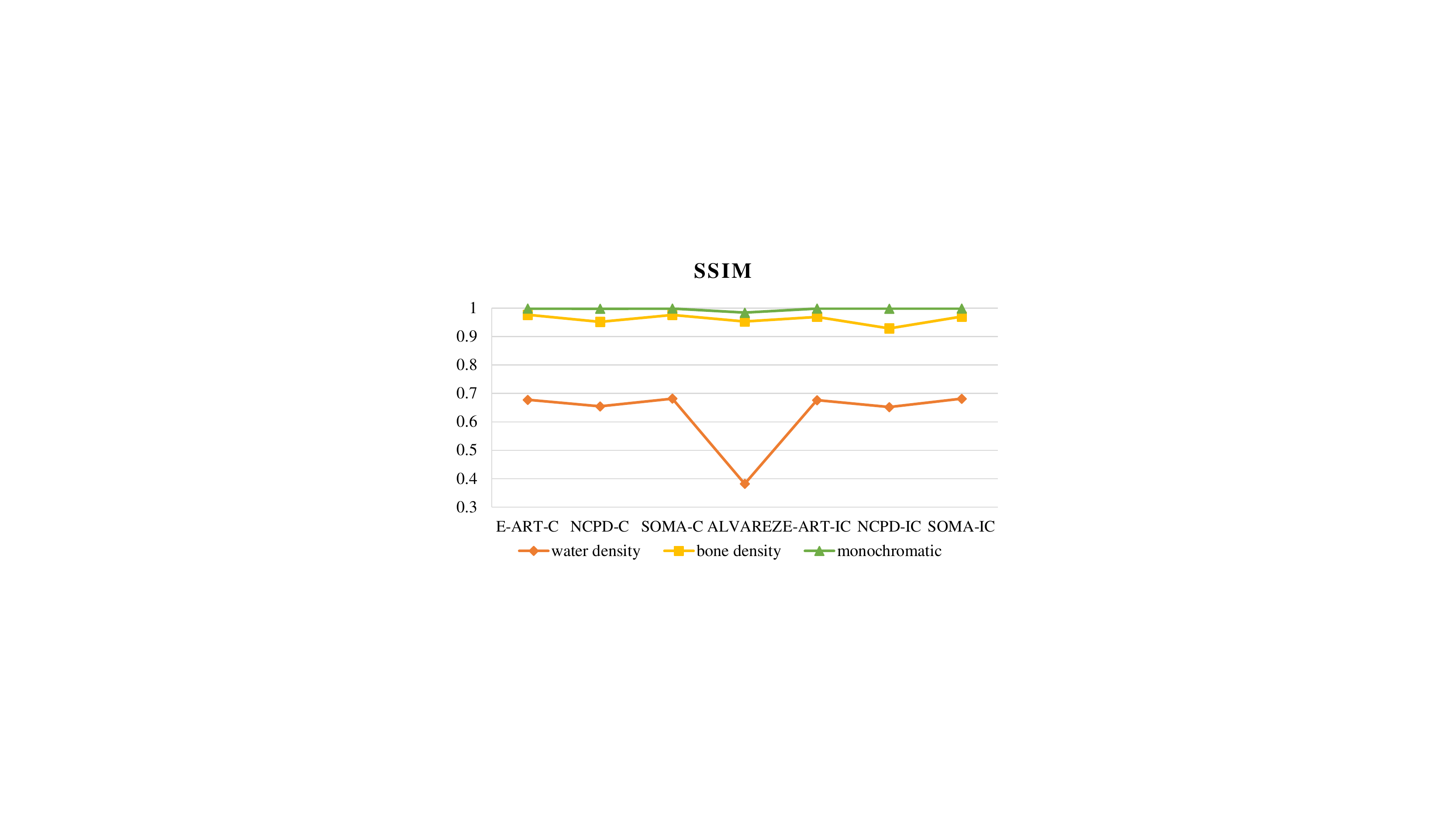}} \hspace{-2cm}
	\vspace{-2.5mm}
	
	\subfigure[]{	
		\label{fig11_c}	
		\includegraphics[width=7.7cm]{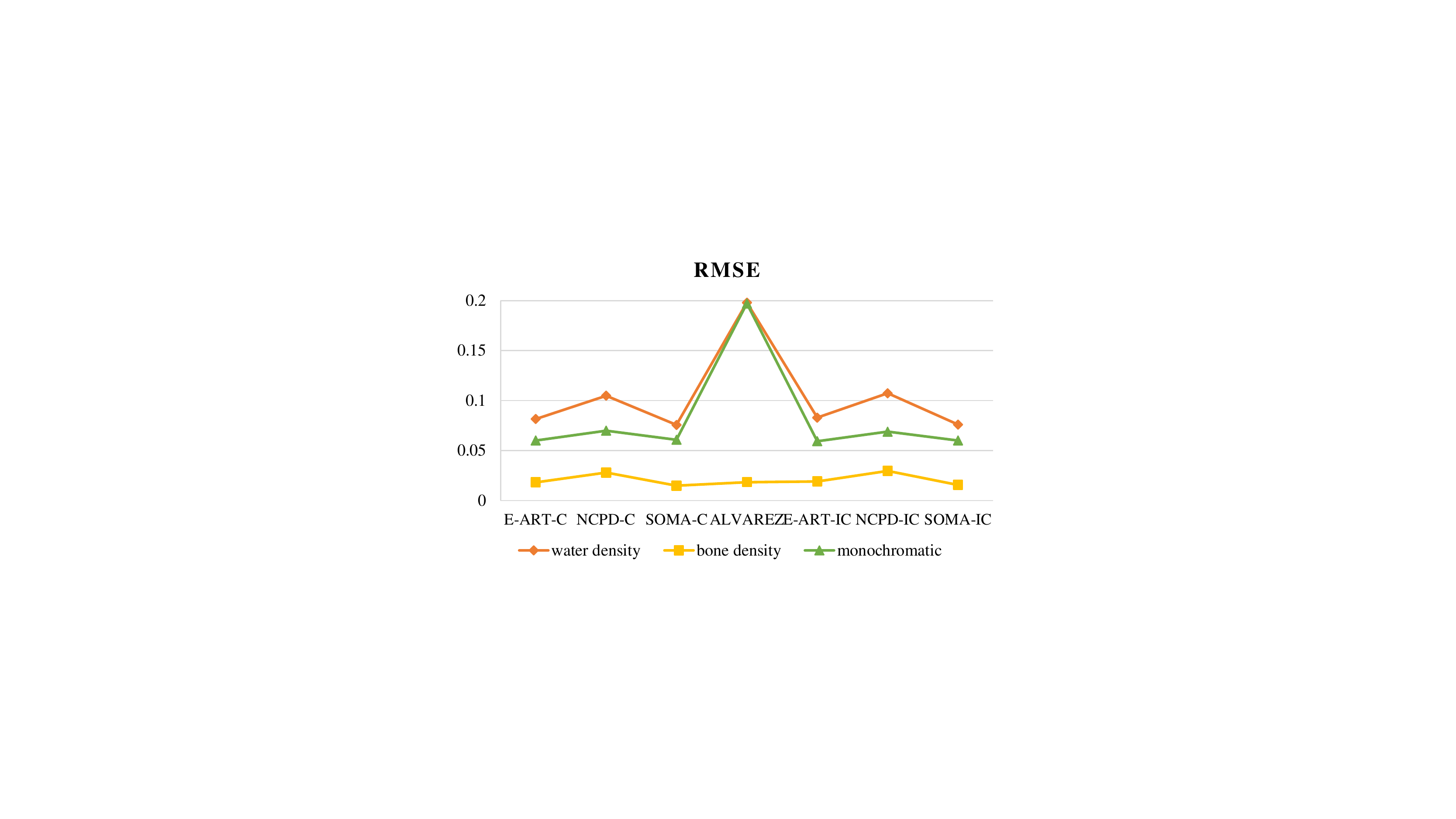}} \hspace{-2cm}
	\vspace{-2.5mm}
	\caption{Quantitative metrics of the noisy data experiment results. (a) PSNR. (b) SSIM. (c) RMSE. X-C represents the results reconstructed from consistent data, X-IC represents from inconsistent data. For the convenience of display, the RMSE of the monochromatic images are enlarged by 50 times.}
	\label{fig11}
\end{figure}

Figure \ref{fig10} shows that, compared with the Alvarez method, the three methods have certain anti-noise and anti-jamming. Besides, after 30 iterations, the three methods can get high-quality bone density images and reconstruct the approximate structures of the water density images. For details in the water density images, the three methods can reconstruct the low contrast areas (marked with blue). However, the spine structures of the water density images (marked with red) reconstructed by the E-ART method and NCPD method are blurred, but the basis material images reconstructed by the proposed method are visually the same as that of the phantom images. Figure \ref{fig11} shows the PSNR, SSIM and RMSE of the reconstruction results in figure \ref{fig10}.

The quantitative metrics of the three methods perform better compared the Alvarez method. The quantitative metrics of the E-ART method are similar to the proposed method. It should be noted that the performance of quantitative metrics is not equivalent to good results. Obviously shown in figure \ref{fig10}, the spine structures of the water density images (marked with red) reconstructed by the E-ART method are different from the phantom. It is reasonable to speculate that some errors in the monochromatic images are offset by the weighted combination of the two density images. That leads to the PSNR of the reconstructed density images of the proposed method is higher than the E-ART method, but the PSNR of the monochromatic images is lower.

Another experiment simulates noisy projections by adding Poisson noise with an initial number of $10^6$ photons to the noise-free data. In the experiment, $D_{\mathrm{data}}^{(n)}$ and $D_{\mathrm{image}}^{(n)}$ of the three methods decrease with the number of iterations increasing. Table 1 list the PSNR and SSIM of the reconstruction results when $D_{\mathrm{image}}^{(n)}<10^{-2}$. Within the range, the quality of reconstruction results is acceptable. 

Observing Table \ref{table1}, there is little difference among the quantitative indicators of the three methods. Further, except for the PSNR of the water density image reconstructed by the NCPD method from geometrically inconsistent data and the SSIM of the bone density image reconstructed by the E-ART method from geometrically consistent data, other quantitative metrics of the results reconstructed by the proposed method are better than the other two methods. On the premise that the reconstruction quality is equivalent, the convergence speed of the proposed method is about 36$\%$-39$\%$ faster than the E-ART method, and 43$\%$-47$\%$ faster than the NCPD method.

\begin{table}[htbp]
	\footnotesize
	\centering
	\caption{SSIM and PSNR of the reconstruction results of noisy data when $D_{\mathrm{image}}^{(n)}<10^{-2}$.}
	\renewcommand\arraystretch{0.7}
	\setlength\tabcolsep{5pt}
	\begin{tabular}{cccccccc}
		\hline
		\specialrule{0em}{1.5pt}{1.5pt}& & \multicolumn{3}{c}{Consistent data} & \multicolumn{3}{c}{Inconsistent data} \\ \cline{3-8}
		\specialrule{0em}{1.5pt}{1pt}& & {\begin{tabular}[c]{@{}c@{}}E-ART 66\\iterations\end{tabular}} & {\begin{tabular}[c]{@{}c@{}}NCPD 74\\ iterations\end{tabular}} & {\begin{tabular}[c]{@{}c@{}}SOMA 39\\ iterations\end{tabular}} & {\begin{tabular}[c]{@{}c@{}}E-ART 74\\ iterations\end{tabular}} & {\begin{tabular}[c]{@{}c@{}}NCPD 83\\ iterations\end{tabular}} & {\begin{tabular}[c]{@{}c@{}}SOMA 47\\ iterations\end{tabular}} \\ \hline
		\specialrule{0em}{1pt}{1pt}\multirow{3}{*}{PSNR} & \begin{tabular}[c]{@{}c@{}}Water density\\image\end{tabular} & 30.499451 & 30.529226 & \textbf{30.548021} & 30.592288 & \textbf{30.648220} & 30.620204 \\
		\specialrule{0em}{1pt}{1pt}& \begin{tabular}[c]{@{}c@{}}Bone density\\image\end{tabular} & 42.139358 & 42.118393 & \textbf{42.155555} & 41.987757 & 41.925252 & \textbf{41.997872} \\
		\specialrule{0em}{1pt}{1pt}& \begin{tabular}[c]{@{}c@{}}Monochromatic\\image\end{tabular} & 69.680017 & 69.421222 & \textbf{69.985931} & 69.889102 & 69.700071 & \textbf{70.052709} \\
		\specialrule{0em}{1pt}{1pt}\multirow{3}{*}{SSIM} & \begin{tabular}[c]{@{}c@{}}Water density\\image\end{tabular} & 0.902462 & 0.899899 & \textbf{0.906882} & 0.904308 & 0.903117 & \textbf{0.906655} \\
		\specialrule{0em}{1pt}{1pt}& \begin{tabular}[c]{@{}c@{}}Bone density\\image\end{tabular} & \textbf{0.997852} & 0.996422 & 0.997780 & 0.997207 & 0.994772 & \textbf{0.997904} \\
		\specialrule{0em}{1pt}{1pt}& \begin{tabular}[c]{@{}c@{}}Monochromatic\\image\end{tabular} & 0.999889 & 0.999879 & \textbf{0.999897} & 0.999892 & 0.999885 & \textbf{0.999896} \\ \hline
	\end{tabular}
	\label{table1}
\end{table}

\subsection{Triple material data experiment}
\label{sect3-3}
\begin{figure}[htbp]
	\centering
	\subfigure[]{	
		\label{fig12_a}
		\includegraphics[width=3.5cm]{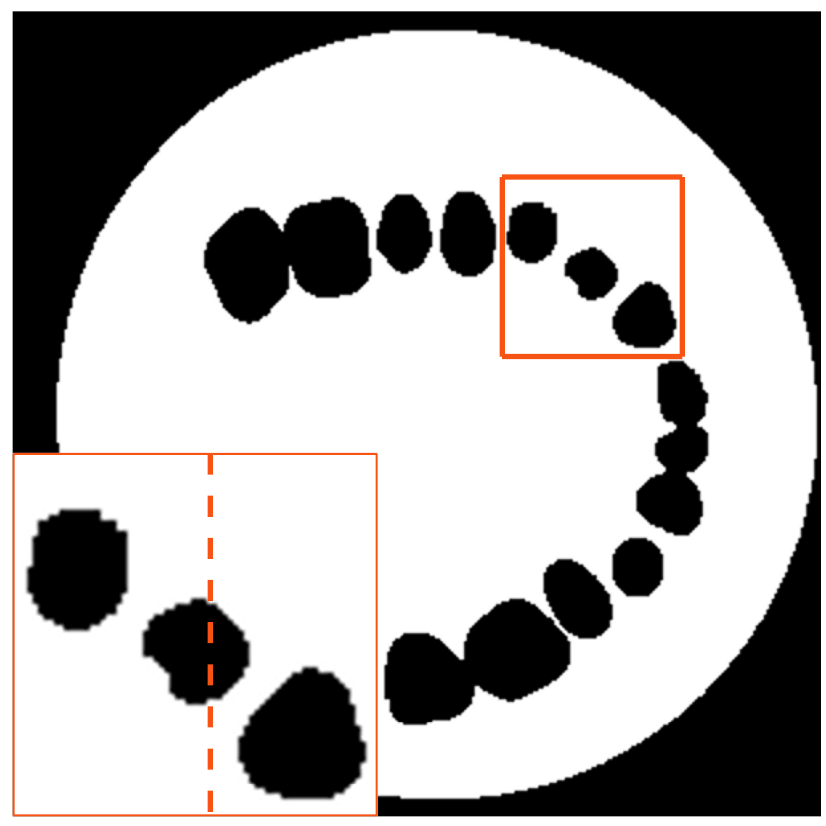}} \hspace{-.2cm}
	\subfigure[]{	
		\label{fig12_b}	
		\includegraphics[width=3.5cm]{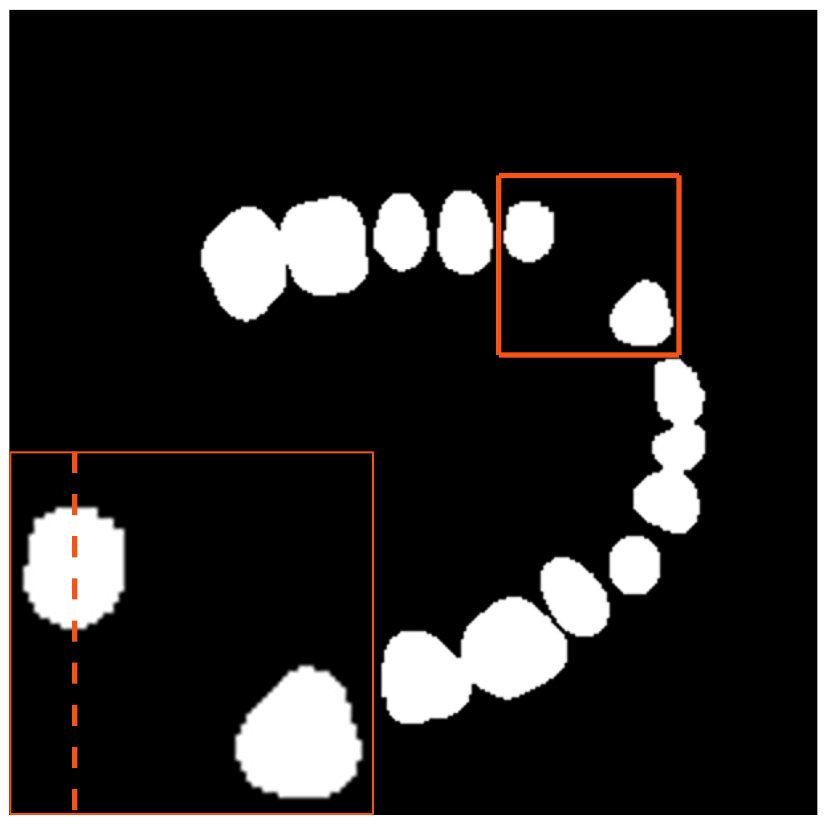}} \hspace{-.2cm}
	\subfigure[]{	
		\label{fig12_c}
		\includegraphics[width=3.5cm]{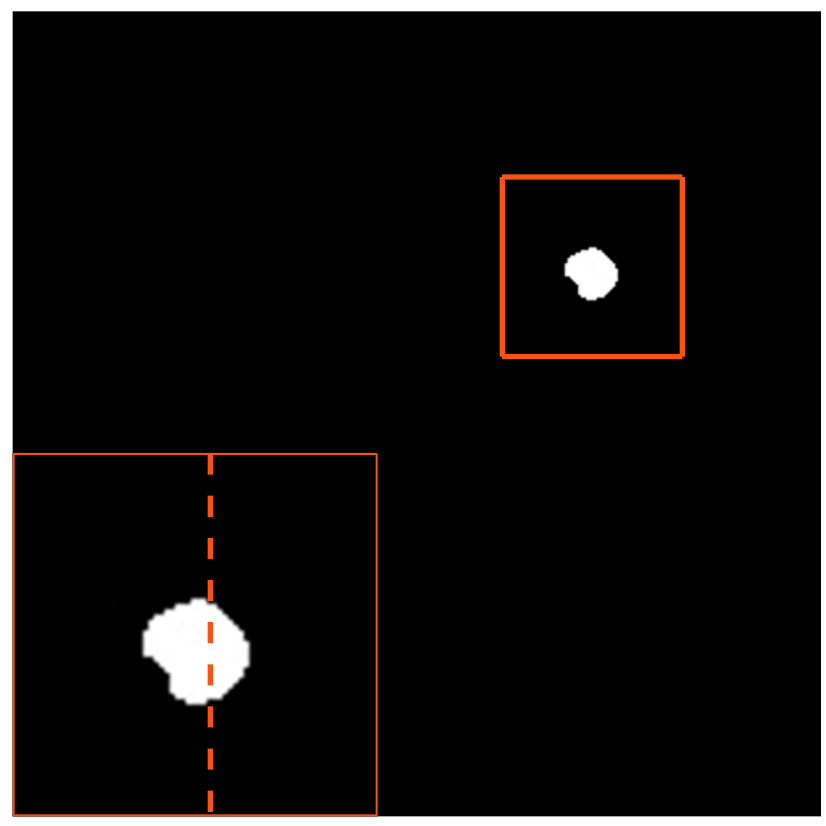}} \hspace{-.2cm}
	\subfigure[]{
		\label{fig12_d}
		\includegraphics[width=3.5cm]{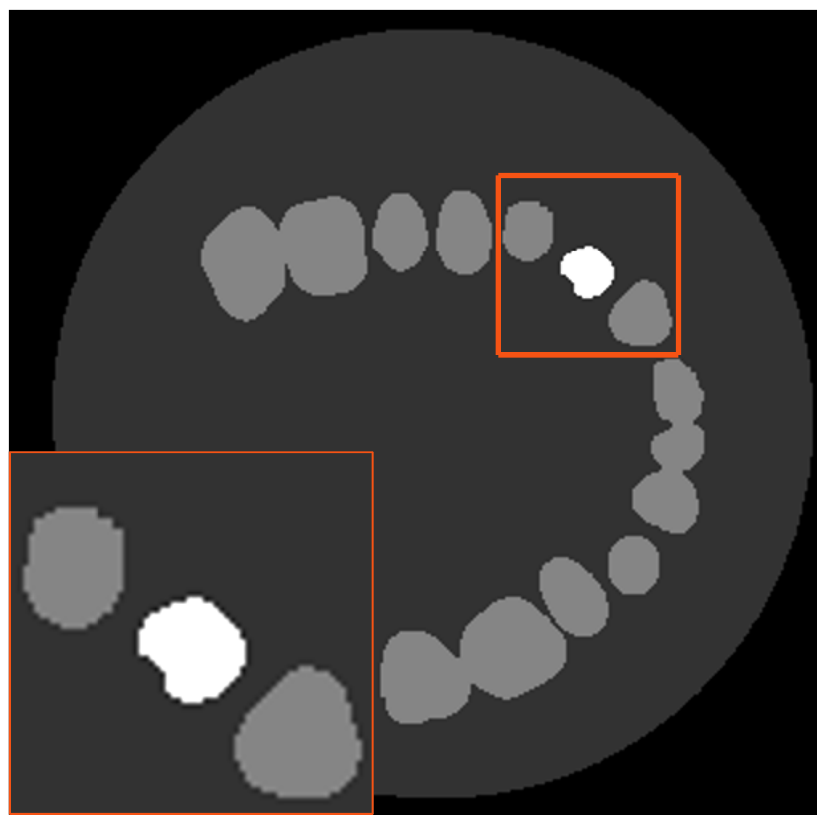}} \hspace{-10mm}
	\vspace{-2.5mm}
	
	\subfigure[]{
		\label{fig12_e}
		\includegraphics[width=7cm]{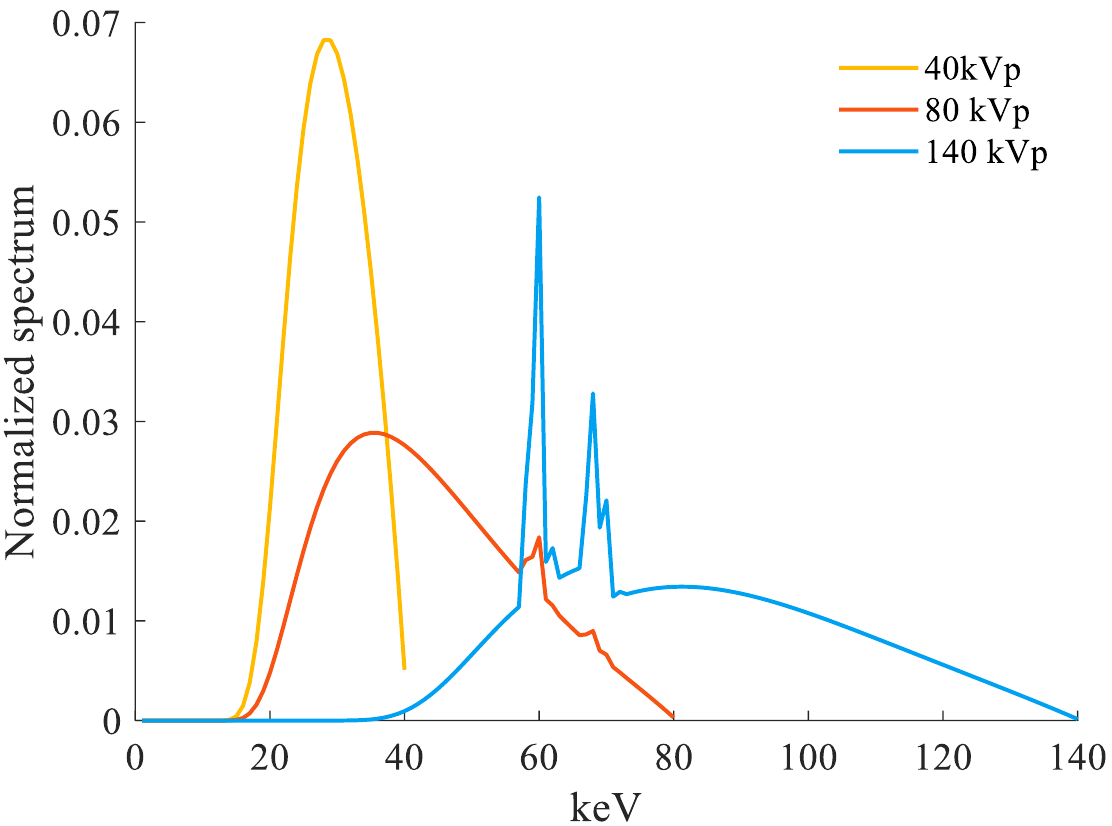}} \hspace{-14mm}
	\vspace{-2.5mm}
	\caption{Phantom and spectra used in triple material data experiment. (a) The water basis material image (Display window: [0, 1.0]). (b) The bone basis material image ([0, 1.0]). (c) The gold basis material image ([0, 1.0]) (d) The monochromatic image at 70keV ([0, 0.1]). (e) Spectra.}
	\label{fig12}
\end{figure}

\begin{figure}[]
	\centering	
	\includegraphics[width=13.8cm]{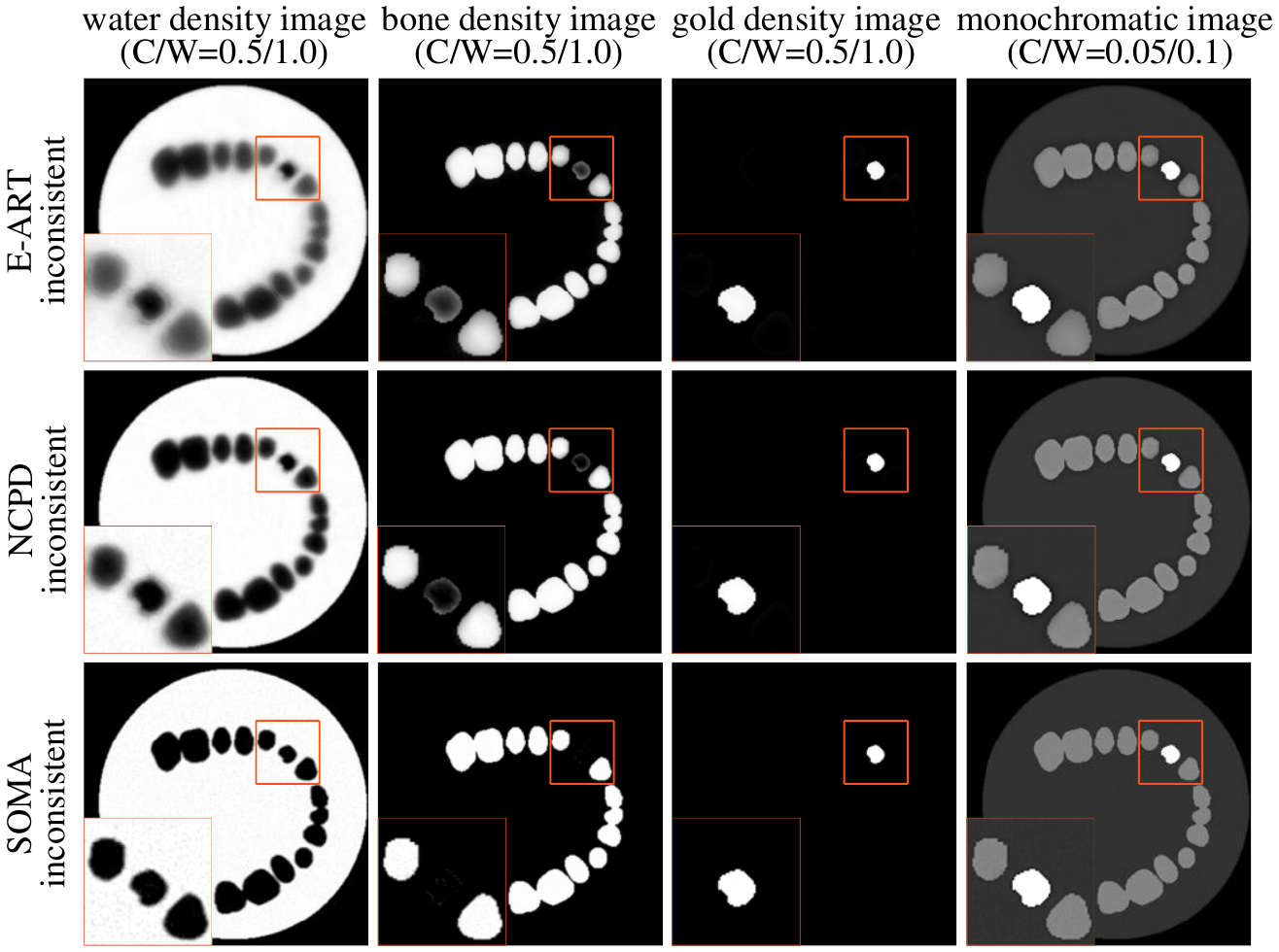} \hspace{-2.4cm}
	\vspace{-2.5 mm}
	\caption{Reconstruction results of triple material data after 200 iterations.}
	\label{fig13}	
\end{figure}

The experiment simulates an oral model included gold teeth with resolution $512\times512$, and the phantom is composed of water, bone, and gold, with densities are 1.0 g/cm$^3$, 1.92 g/cm$^3$, and 19.32 g/cm$^3$ respectively. Figure \ref{fig12}(a)-\ref{fig12}(d) shows the basis material images and the monochromatic image at 70 keV energy. Using Spectrum GUI to simulate the spectra, the X-ray source is GE Maxiray 125 X-ray tube, the tube voltage is set to 40 kVp, 80 kVp, and 140 kVp, a 1 mm copper filter is added in front of the X-ray source at 140 kVp, and the normalized X-ray spectra are shown in figure \ref{fig12}(e). The mass attenuation coefficients of water, bone, and gold can be obtained from the NIST website. The scan configurations are the same as those in section \ref{sect3-1}. Using the standard full-scan configuration, and a total of 720 projections are collected under each spectrum. When simulating the process of generating geometrically inconsistent data, the initial angle of the single full scanning differs by 0.16$^{\circ}$. Adding Poisson noise with an initial number of $10^6$ photons to the geometrically inconsistent data. Figure \ref{fig13} shows the reconstruction results of the E-ART method and the proposed method after 200 iterations.

After 200 iterations, the three methods can accurate the gold density images. The structure of the water density images and bone density images, which reconstructed with the E-ART method and the NCPD method, has errors and these errors are offset by the weighted combination, therefore, there is no visual difference between the monochrome images reconstructed by the three methods. It takes thousands of iterations for the E-ART method and the NCPD method to get the correct basis material images, which is because of the greater overlap between the spectra used in triple material data experiments and the stronger ill-condition of the nonlinear model. The proposed method uses Schmidt orthogonalization to find the optimal solution along the orthogonal direction orthogonal, which helps to obtain high-precision reconstruction results with fewer iterations.

\begin{figure}[htbp]
	\centering
	\subfigure{
	\begin{minipage}[t]{0.49\linewidth}
		\centering
		\centerline{\includegraphics[height=4.5cm]{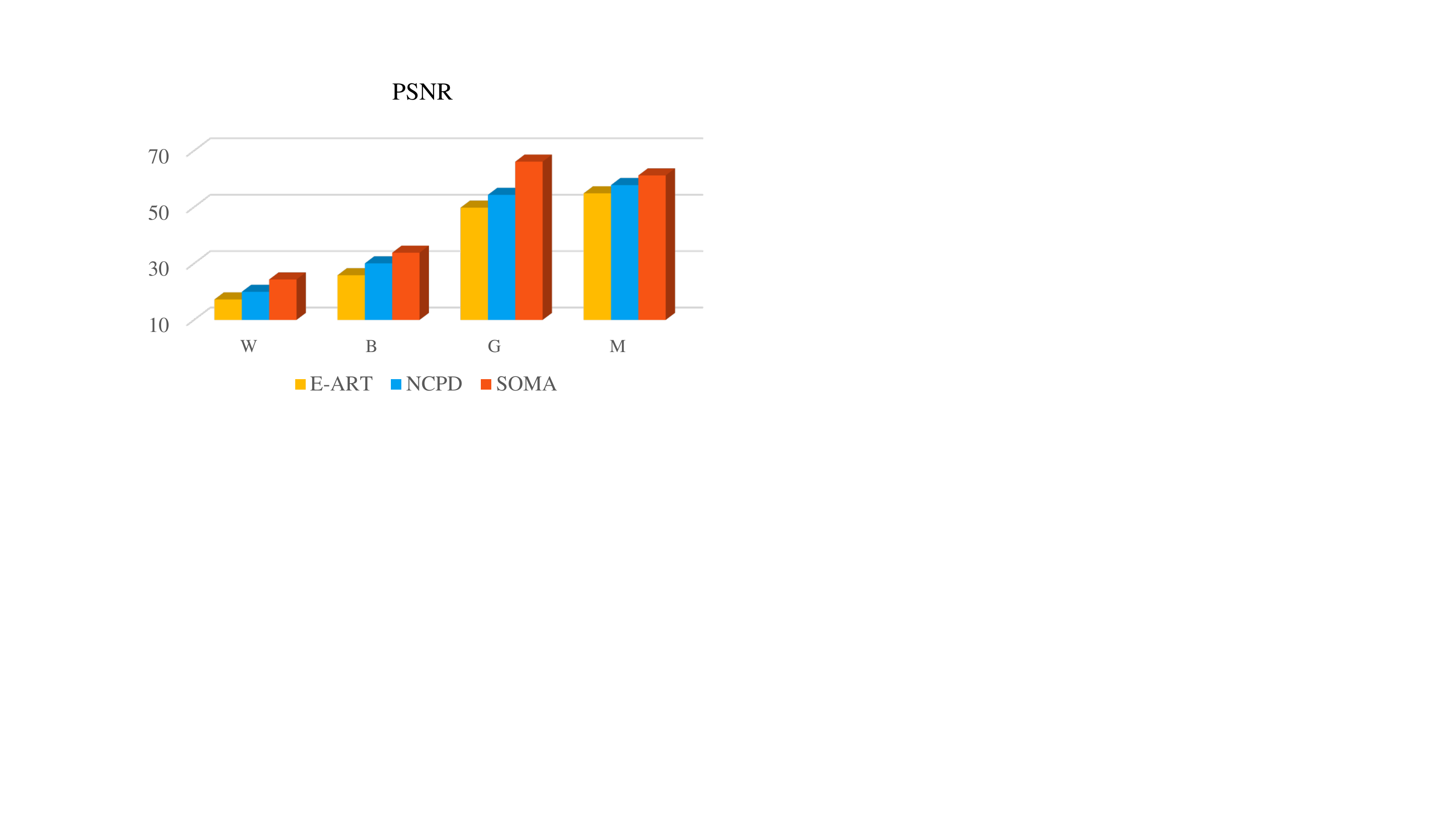}\label{fig14_a}}
		\centerline{(a)}
		\vspace{0.25cm}
		\centerline{\includegraphics[height=4.5cm]{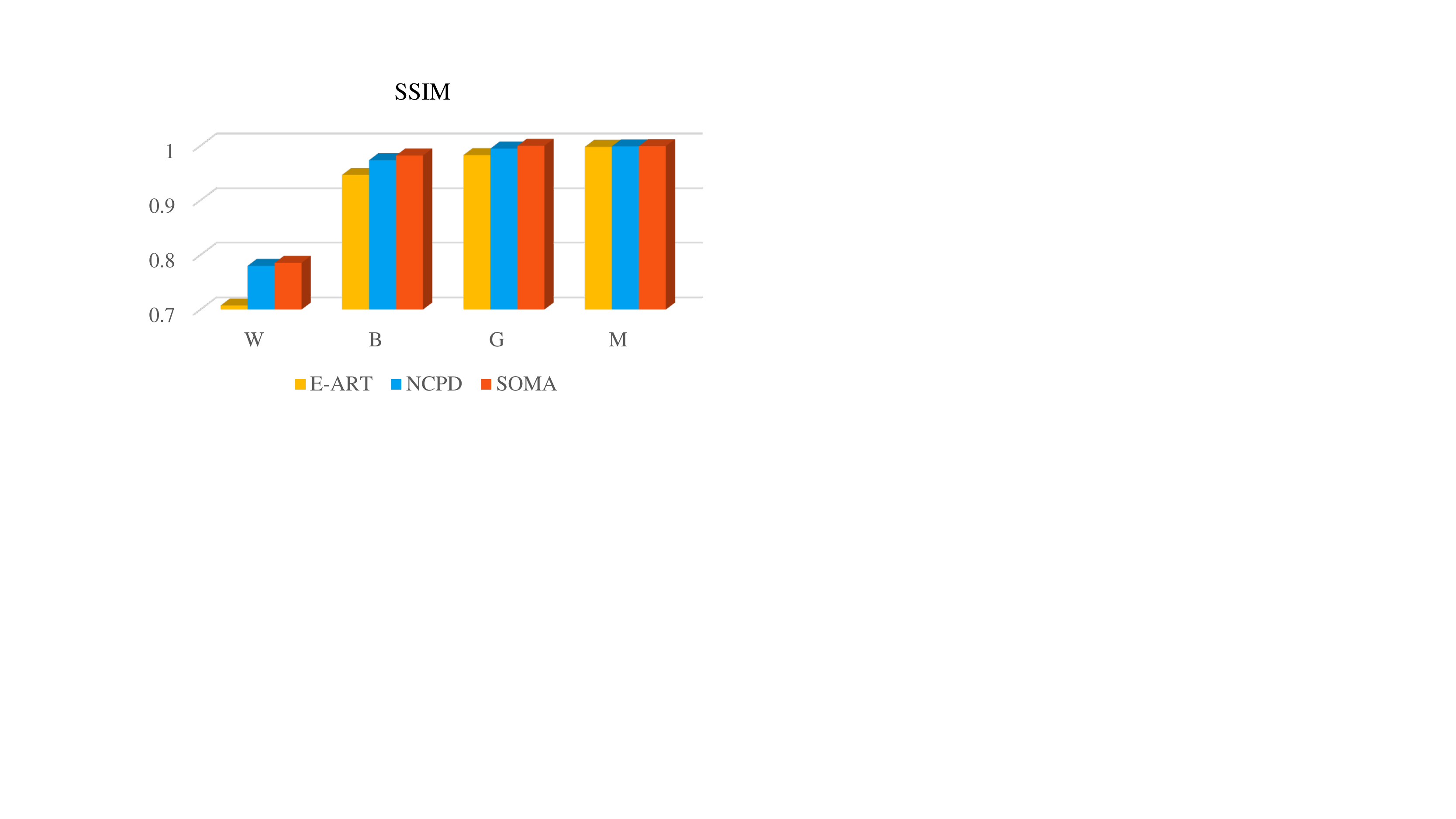}\label{fig14_b}}
		\centerline{(b)}
		\vspace{0.25cm}
		\centerline{\includegraphics[height=4.5cm]{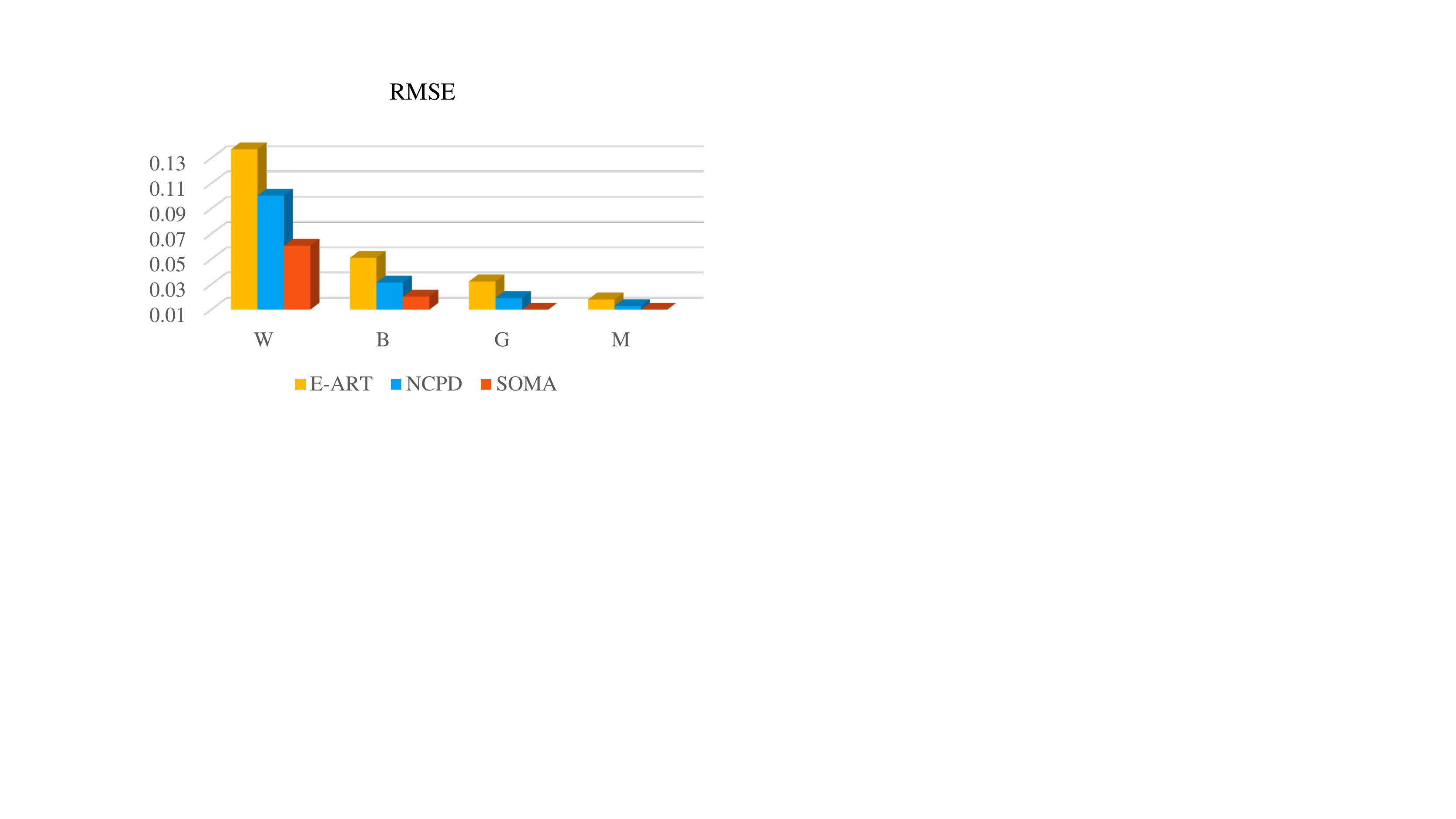}\label{fig14_c}}
		\centerline{(c)}
		\vspace{0.25cm}
	\end{minipage}} \hspace{-2mm}
	\subfigure{
	\begin{minipage}[t]{0.49\linewidth}
		\centering
		\centerline{\includegraphics[height=4.5cm]{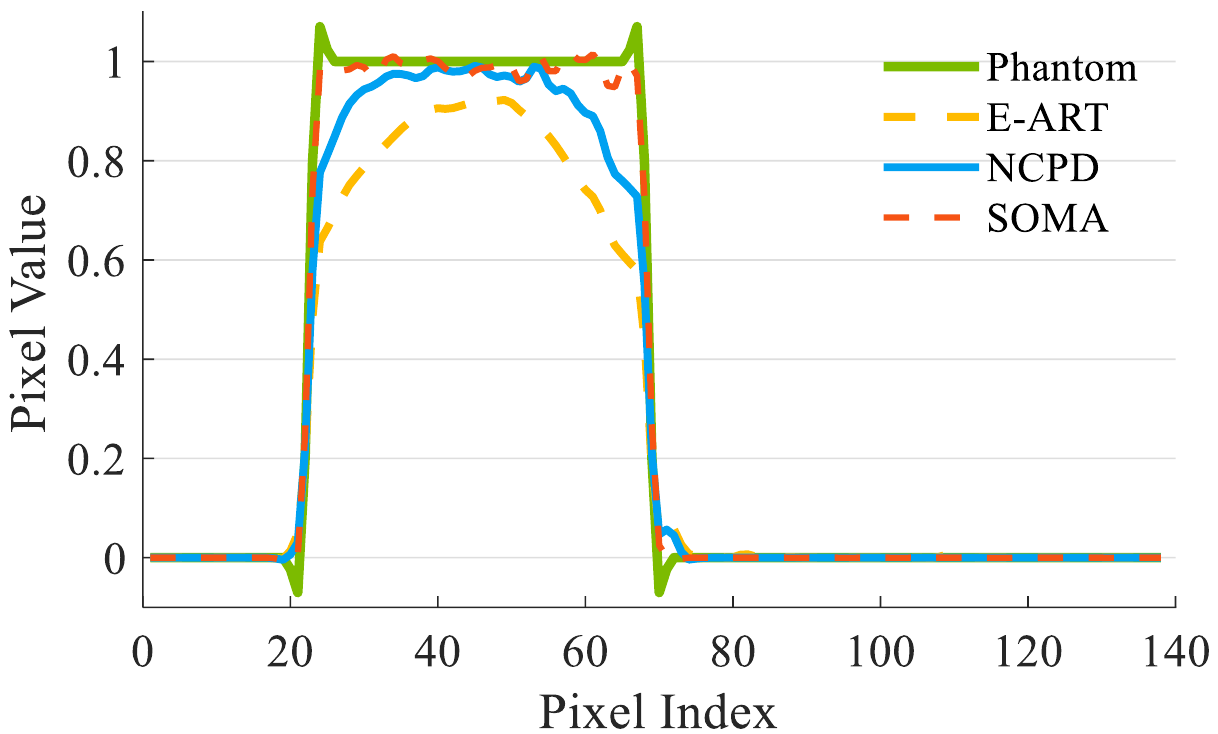}\label{fig14_d}}
		\centerline{(d)}
		\vspace{0.25cm}
		\centerline{\includegraphics[height=4.5cm]{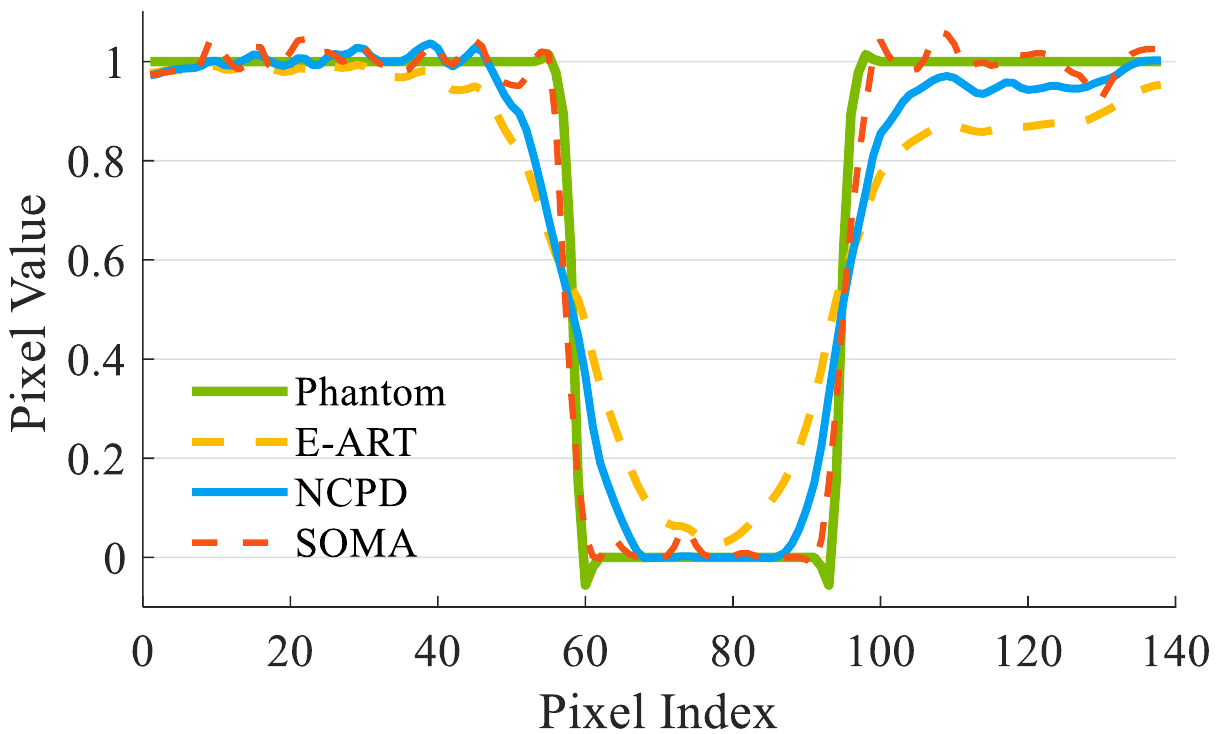}\label{fig14_e}}
		\centerline{(e)}
		\vspace{0.25cm}
		\centerline{\includegraphics[height=4.5cm]{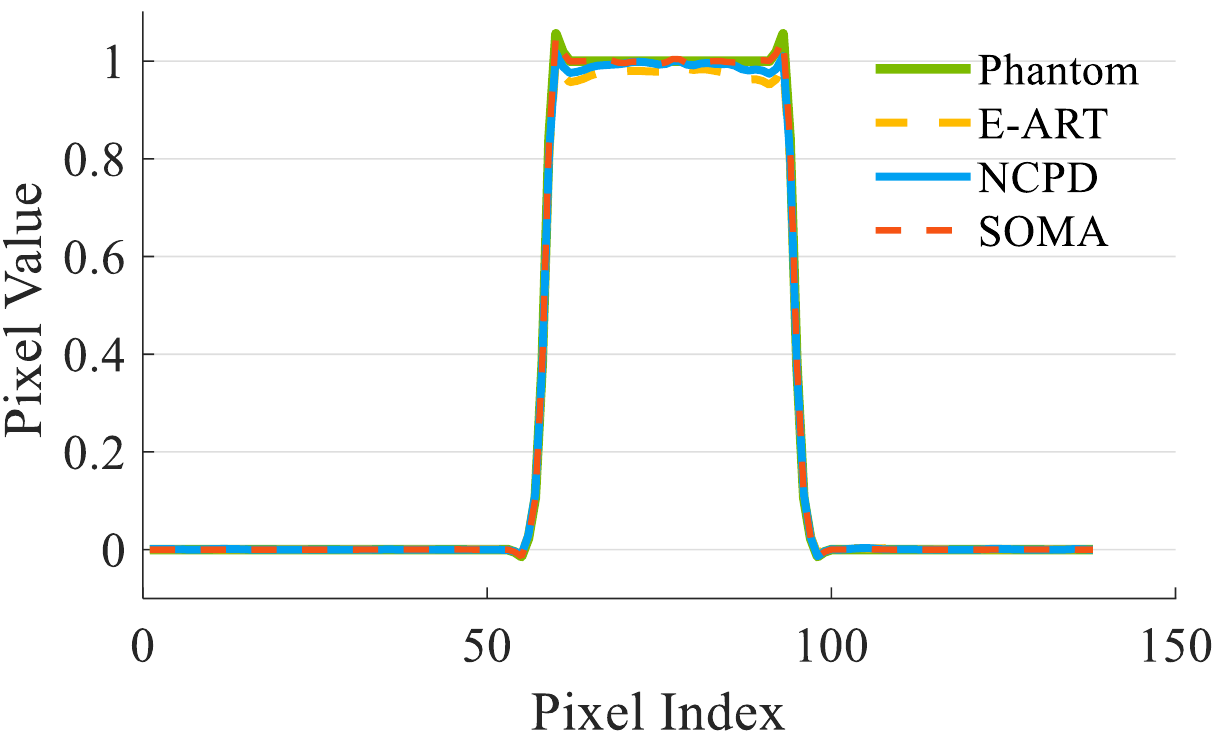}\label{fig14_f}}
		\centerline{(f)}
		\vspace{0.25cm}
	\end{minipage}}\hspace{-1.8cm}
	\vspace{-5mm}
	\caption{Quantitative metrics and profiles of the triple material data experiment results. (a) PSNR. (b) SSIM. (c) RMSE. (d) Profiles of the water density image (marked with red line in figure \ref{fig12}(a)). (e) Profiles of the bone density image (marked with red line in figure \ref{fig12}(b)). (f) Profiles of the gold density image (marked with red line in figure \ref{fig12}(c)). W represents the reconstructed water density image, B represents the bone density image, G represents the gold density image and M represents the monochromatic images at 70keV. For the convenience of the display, the RMSE of the gold density image and the monochromatic images are enlarged by 10 times.}
	\label{fig14}	
\end{figure}

Figure \ref{fig14} shows some quantitative metrics and profiles of the reconstruction results. Observing figure \ref{fig14}(a)-\ref{fig14}(c), the quantitative indicators and it can be observed that the PSNR and SSIM of the proposed method are higher than the E-ART method and the NCPD method, while the RMSE is lower than two methods than the other two methods. Observing figure \ref{fig14}(d)-\ref{fig14}(f), the profiles of the reconstructed results show that the value of the basis material images reconstructed by the proposed method is closer to the phantom.

Continue the iterations of the E-ART method and the NCPD method until $D_{\mathrm{image}}^{(n)}<10^{-2}$, the E-ART method stops the iteration at 1250, the NCPD method stops at 1036 and the proposed method stops at 231. In this case, when the reconstruction quality is equivalent, the convergence speed of the proposed method is about 81$\%$ faster than the E-ART method, and about 78$\%$ faster than the NCPD method.

\subsection{Real data experiment}
\label{sect3-4}
\begin{figure}[htbp]
	\centering
	\subfigure[]{	
		\label{fig15_a}
		\includegraphics[width=7cm]{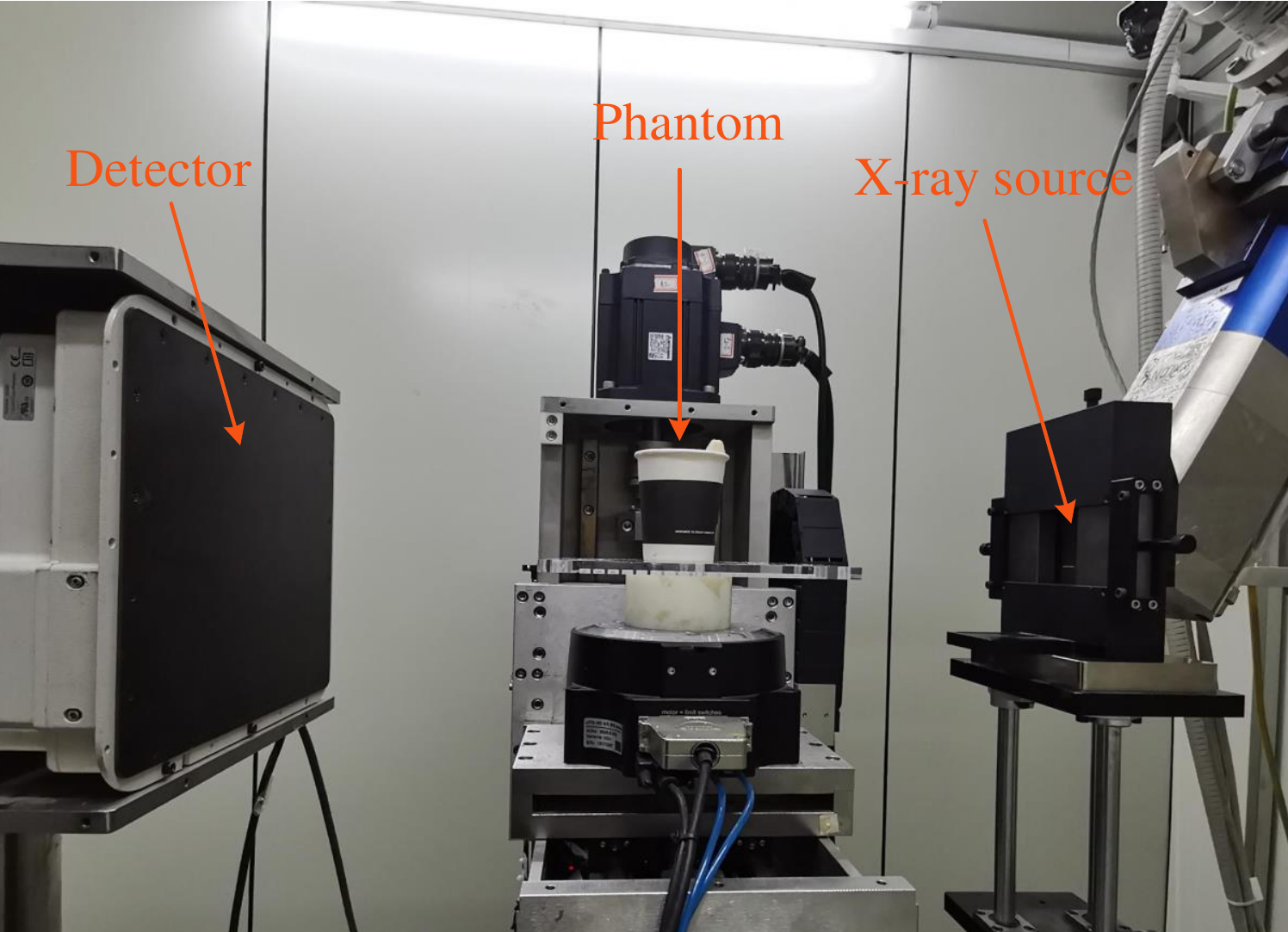}} \hspace{5mm}
	\subfigure[]{	
		\label{fig15_b}	
		\raisebox{0.1\height}{\includegraphics[width=4cm]{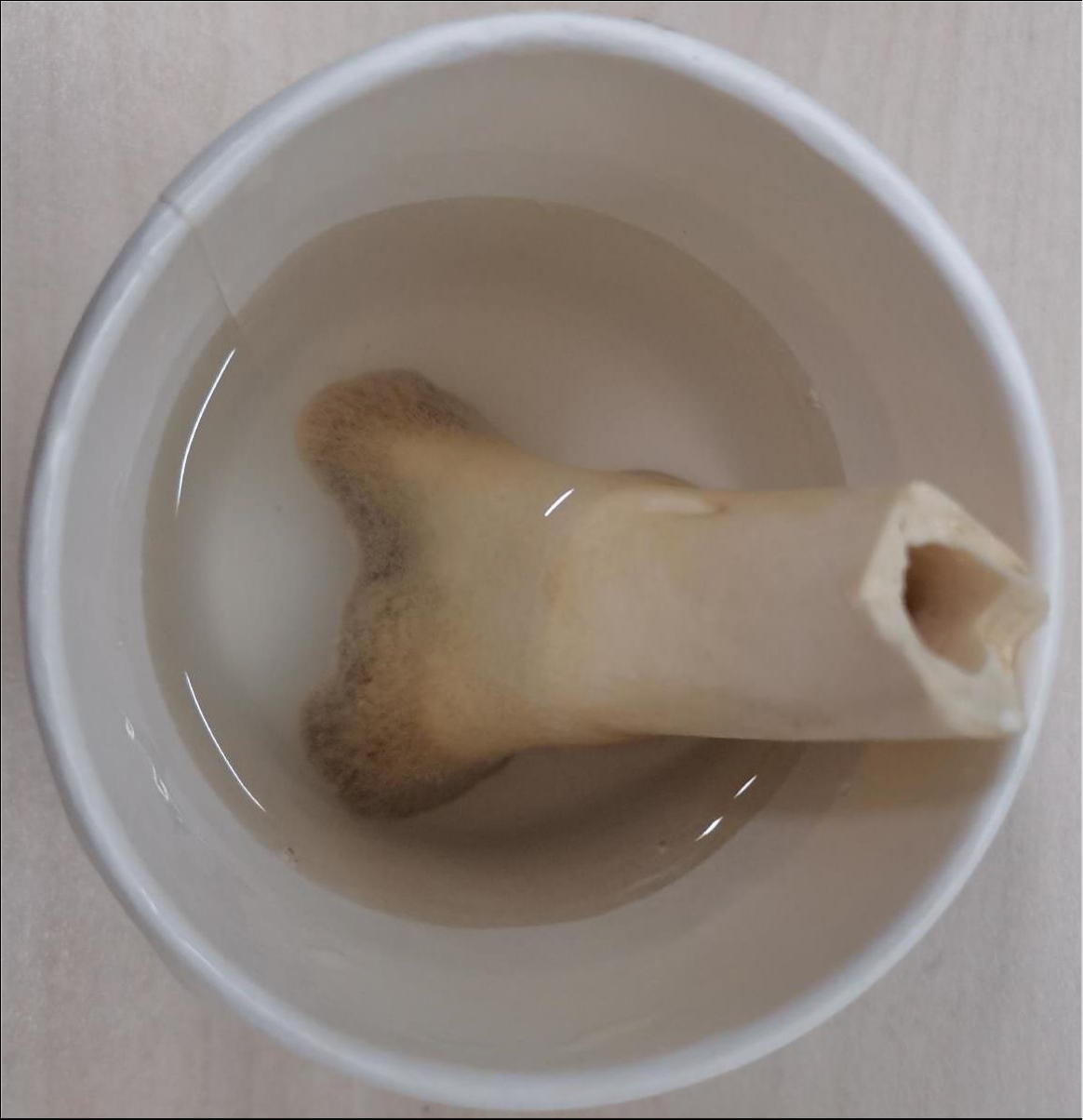}}} \hspace{-2cm}
	\vspace{-2.5mm}
	
	\subfigure[]{	
		\label{fig15_c}
		\includegraphics[width=7cm]{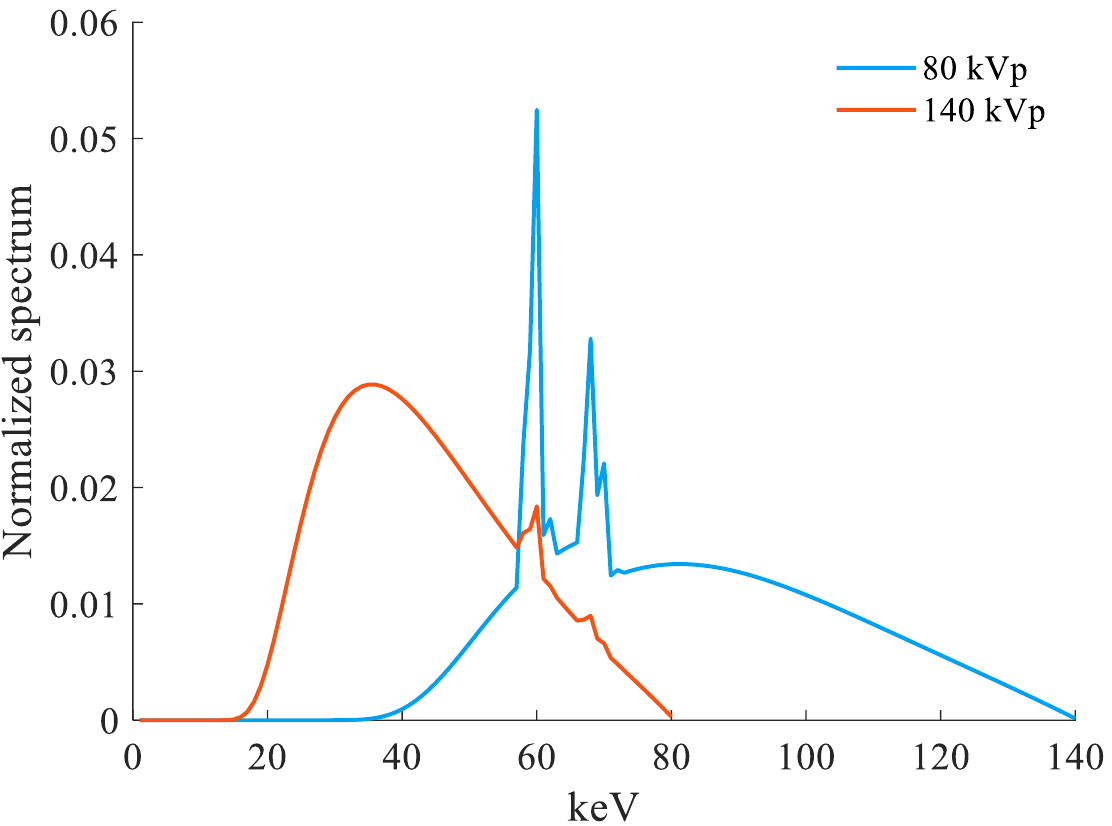}} \hspace{-2cm}
	\vspace{-2.5mm}
	\caption{CT system, phantom, and spectra for the real data experiment. (a) Photograph of the industrial CT system in our laboratory. (b) Bone-water phantom. (c) The estimated spectra.}
	\label{fig15}
\end{figure}

In this section, the proposed method is performed on real data with complex structures. The equipment, phantom, and spectra are shown in figure \ref{fig15}(a)-\ref{fig15}(c), the mass attenuation coefficients of water and bone can be obtained from the NIST website, and the scan configurations are shown in Table \ref{table2}.

\begin{table}[htbp]
	\footnotesize
	\centering
	\caption{Scan configuration for the real data experiment.}
	\renewcommand\arraystretch{1.1}
	\setlength\tabcolsep{10pt}
	\begin{tabular}{ccc}
		\hline
		\specialrule{0em}{1pt}{1pt}& Scan 1 & Scan 2 \\
		\hline
		Voltage & 80 kVp & 140 kVp \\
		Current & 240 uA & 120 uA \\
		Filter & 1.5 mm Al & 0.5 mm Cu \\
		Exposure time per projection & 0.5 s & 0.6 s \\
		SOD & 355.61 mm & 355.61 mm \\
		SDD & 673.96 mm & 673.96 mm \\
		Projections & 1440 & 1440 \\
		\hline
	\end{tabular}
	\label{table2}
\end{table}

\begin{figure}[]
	\centering	
	\includegraphics[width=13.8cm]{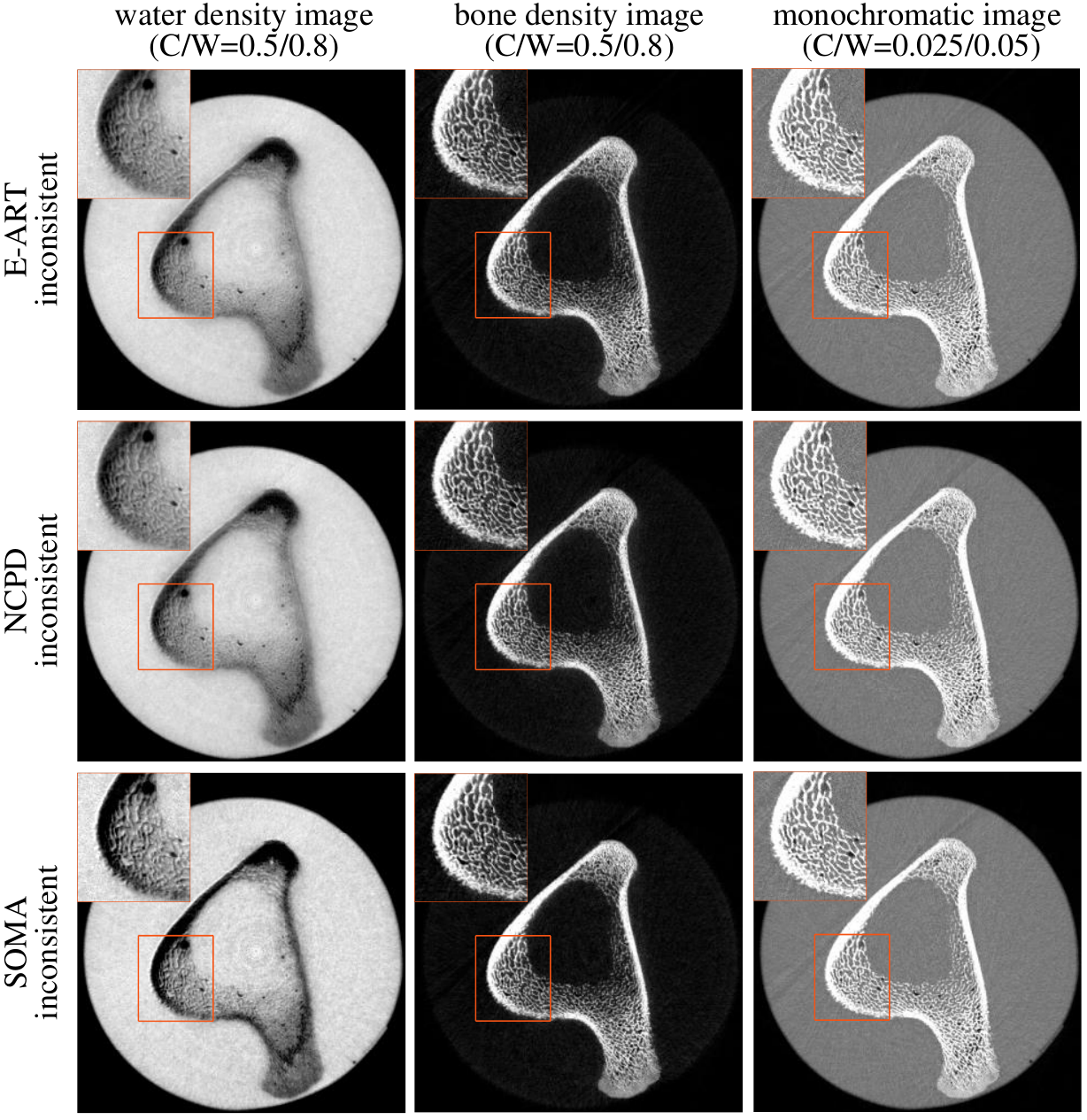} \hspace{-2cm}
	\vspace{-2.5 mm}
	\caption{Reconstruction results of real data after 5 iterations.}
	\label{fig16}	
\end{figure}

This experiment takes the iterations as the stopping criterion, and the maximum iteration is set to 5. The parameter of the proposed method are as $\beta=0.05$, $\epsilon=10^{-8}$, $\eta = 1.5$, $\kappa = 0.95$ and $\alpha = 1.0$. Figure \ref{fig16} shows the reconstruction results of the E-ART method, the NCPD method and the proposed method when the stopping criterion is met.

Figure \ref{fig16} shows that when the iterations are the same, the three methods can get acceptable the bone density images. The bone trabecular structures (marked with red) in the water basis material images reconstructed by the E-ART method and the NCPD method are unclear, and on the contrary, the structures are clearly visible in the results reconstructed by the proposed method. The upper left part of the bone trabecula in the water density image has an obvious edge, and the other two methods will get similar results if they continue to iterate. The other two methods and the proposed methods can deal with geometrically inconsistent data, thus the edge is not caused by solution error. It is speculated empirically to be caused by the movement of the object during the scanning process. The results illustrate that the proposed method can process real data and has practical value.

\section{Discussion}
\label{sect4}
This paper proposed a general iterative method, the so-called SOMA method, to invert the nonlinear equations. The core of the proposed method is getting the orthogonal direction by Schmidt orthogonalization, which accelerates the convergence speed greatly.

In section \ref{sect3}, the validity of the proposed method is verified by MSCT basis material decomposition experiments. Three simulation data experiments verify the numerical convergence of the proposed method, the robustness to noise, and the feasibility of the multiple basis material decomposition. A real data experiment illustrates the practical value of the proposed method. The above experiment results show that the convergence speed of the proposed method is faster than the E-ART method and the NCPD method with high-precision solutions.

The proposed method is not sensitive to most of the parameters. In order to make parameter selection easier, an adaptive step size strategy is given in section \ref{sect2-3}. It should be emphasized that the adaptive step size strategy is not necessary. Fixed step size can satisfy the need when the reconstruction quality requirement is not high. Besides the adaptive step size strategy, $\beta $ can be attenuated as the iterations increase, one of the attenuation strategies is exponential, the corresponding formula is:
\begin{equation}
	\label{eq-beta}
	\beta = \beta_0\cdot\kappa^{\frac{n-1}{N}},
\end{equation}
where $\beta_0 $ is the initial value, $\kappa(\in(0, 1])$ is the attenuation ratio, $n$ is the current iteration, $N$ is the total iterations. The adaptive step size strategy and (\ref{eq-beta}) are for reference only, other strategies such as linear method or piecewise linear method can also be used.

In practical applications, the optimization problem combined with the assumption of minimizing the total variation can be proposed. We have tested a variety of optimization models, and the proposed method combines well with the ASD-NC-POCS method \cite{ref26}. Refer to the ASD-NC-POCS method, an optimization model can be constructed as
\begin{equation}
	\label{eq-opmodel}
	\begin{array}{l}
		f^{*} = \arg \min \sum_{m=1}^{M} \Vert f_m \Vert_{TV} \quad s.t. \quad \sqrt{\frac{\sum_{k=1}^{K}\Vert \tilde{p}_k - p_k\Vert^2_2}{\sum_{k=1}^{K}\Vert \tilde{p}_k \Vert^2_2}}<\tau,
	\end{array}
\end{equation}
where $\tau$ is the threshold of the data term and $\tilde{p}_k$ is the measured polychromatic projection of the $k$-th spectrum. For the above optimization model, we give a natural but non-strict solving method. First, use the proposed method to solve the data term (or non-linear model), then perform the POCS method, and at last solve the total variation minimization by the steepest descent method. (\ref{eq-opmodel}) is not the only model that can be combined with the proposed method and the given method is not the only method to solve (\ref{eq-opmodel}). Other optimization models and solutions are worthy of in-depth study but not discussed further here.

In this paper, the influence of scattered photons is ignored when modeling MSCT reconstruction problem, however, the scattered photons have a great influence in the actual scanning. MSCT model containing scattered photons and its solution is the focus of the next step.

\section*{Acknowledgments}
This work was supported by National Natural Science Foundation of China (No.61827809) and the National Key Research and Development Program of China (No.2020YFA0712200)

\appendices
\section*{Appendix A.}
\label{Appendix A}
\setcounter{section}{1}
\setcounter{table}{0}
\setcounter{figure}{0}
\setcounter{equation}{0}
\renewcommand\theequation{A.\arabic{equation}}
\renewcommand{\thetable}{A.\arabic{table}}
\renewcommand{\thefigure}{A.\arabic{figure}}
In this section, the recursion process of the Schmidt orthogonalization is shown. The initial vectors are ${\bf g}_1, {\bf g}_2, \cdots, {\bf g}_K$ and the orthogonal vectors are ${\bf d}_1, {\bf d}_2, \cdots, {\bf d}_K$.

Firstly, let ${\bf d}_1 = {\bf g}_1 = {\bf P}_1{\bf g}_1$ and ${\bf P}_1$ is set as the unit matrix. Then, the orthogonal vectors ${\bf d}_2$ can be obtained by the Schmidt orthogonalization
\begin{equation}
	\label{eq-so-d2}
	\begin{aligned}
		{\bf d}_2 &= {\bf g}_2 - {\bf d}_1\frac{<{\bf g}_2,{\bf d}_1>}{<{\bf d}_1,{\bf d}_1>} \\
		&= {\bf g}_2 - {\bf d}_1\frac{{\bf d}_1^\top {\bf g}_2}{{\bf d}_1^\top {\bf d}_1} \\
		&= ({\bf I}-\frac{{\bf d}_1{\bf d}_1^\top}{{\bf d}_1^\top {\bf d}_1}){\bf g}_2 \\
		&= ({\bf P}_1-\frac{{\bf d}_1{\bf d}_1^\top}{{\bf d}_1^\top {\bf d}_1}){\bf g}_2 \\
		&= {\bf P}_2{\bf g}_2,
	\end{aligned}
\end{equation}
where $<{\bf a},{\bf b}>$ represents the dot product of ${\bf a}$ and ${\bf b}$.

Next, the third orthogonal vectors ${\bf d}_2$ can be calculated in the same way
\begin{equation}
	\label{eq-so-d3}
	\begin{aligned}
		{\bf d}_3 &= {\bf g}_3 - {\bf d}_1\frac{<{\bf g}_3,{\bf d}_1>}{<{\bf d}_1,{\bf d}_1>} -  {\bf d}_2\frac{<{\bf g}_3,{\bf d}_2>}{<{\bf d}_2,{\bf d}_2>}\\
		&= {\bf g}_3 - {\bf d}_1\frac{{\bf d}_1^\top {\bf g}_3}{{\bf d}_1^\top {\bf d}_1} - {\bf d}_2\frac{{\bf d}_2^\top {\bf g}_3}{{\bf d}_2^\top {\bf d}_2} \\
		&= ({\bf I}-\frac{{\bf d}_1{\bf d}_1^\top}{{\bf d}_1^\top {\bf d}_1} -\frac{{\bf d}_2{\bf d}_2^\top}{{\bf d}_2^\top {\bf d}_2}){\bf g}_3 \\
		&= ({\bf P}_2 -\frac{{\bf d}_2{\bf d}_2^\top}{{\bf d}_2^\top {\bf d}_2}){\bf g}_3 \\
		&= {\bf P}_3{\bf g}_3,
	\end{aligned}
\end{equation}

Similarly, it can be obtained by recursion
\begin{equation}
	\label{eq-so-dk}
	\begin{aligned}
		{\bf d}_k &= {\bf g}_k - {\bf d}_1\frac{<{\bf g}_k,{\bf d}_1>}{<{\bf d}_1,{\bf d}_1>} - \cdots - {\bf d}_{k-1}\frac{<{\bf g}_{k},{\bf d}_{k-1}>}{<{\bf d}_{k-1},{\bf d}_{k-1}>}\\
		&= ({\bf I}-\frac{{\bf d}_1{\bf d}_1^\top}{{\bf d}_1^\top {\bf d}_1} -\cdots-\frac{{\bf d}_{k-1}{\bf d}_{k-1}^\top}{{\bf d}_{k-1}^\top {\bf d}_{k-1}}){\bf g}_k \\
		&= ({\bf P}_{k-1} -\frac{{\bf d}_{k-1}{\bf d}_{k-1}^\top}{{\bf d}_{k-1}^\top {\bf d}_{k-1}}){\bf g}_k \\
		&= {\bf P}_k{\bf g}_k.
	\end{aligned}
\end{equation}

So, the formula to get the orthogonal direction is
\begin{equation}
	\label{eq-so-dk-final}
	{\bf d}_k = {\bf P}_{k}{\bf g}_k,
\end{equation}
and the update process of the orthogonal correction matrix is 
\begin{equation}
	\label{eq-so-pk}
	{\bf P}_k = {\bf P}_{k-1}-\frac{{\bf d}_{k-1}{\bf d}_{k-1}^\top}{{\bf d}_{k-1}^\top {\bf d}_{k-1}}.
\end{equation}
In order to avoid the zero denominator, a very small value $\epsilon$ is added in the denominator as the correction term. As a result, the final update formula is
\begin{equation}
	\label{eq-so-pk-mod}
	{\bf P}_k = {\bf P}_{k-1}-\frac{{\bf d}_{k-1}{\bf d}_{k-1}^\top}{{\bf d}_{k-1}^\top {\bf d}_{k-1}+\epsilon}.
\end{equation}

\section*{Appendix B.}
\label{Appendix B}
\setcounter{section}{2}
\setcounter{table}{0}
\setcounter{figure}{0}
\setcounter{equation}{0}
\renewcommand\theequation{B.\arabic{equation}}
\renewcommand\thetable{B.\arabic{table}}
\renewcommand\thefigure{B.\arabic{figure}}

This section gives a simple brief convergence proof of the SOMA method. 

The convergence proof is equivalent to proving that the sequence $\{{\bf x}_k\}$ generated by ${\bf x}_k = {\bf x}_{k-1}+\alpha_k{\bf d}_{k}$ converges to the solution ${\bf x}^*$ of the linear system ${\bf A}{\bf x}={\bf b}$.

There are a series of linearly independent vectors ${\bf d}_1,{\bf d}_2,\cdots,{\bf d}_K$ and they span the solution space. The difference between ${\bf x}^*$ and ${\bf x}_0$ can be written as
\begin{equation}
	\label{eq-B-x*-x}
	{\bf x}^* - {\bf x}_0 = a_1{\bf d}_1 + a_2{\bf d}_2 + \cdots + a_K{\bf d}_K,
\end{equation}
where $a_k$ is a scalar. Multiply both sides of the formula by ${\bf d}_k^\top {\bf A}$ and use the property of orthogonality, then
\begin{equation*}
	\label{eq-B-dkA(x*-x)}
	\begin{aligned}
		{\bf d}_k^\top {\bf A}({\bf x}^* - {\bf x}_0) &= {\bf d}_k^\top {\bf A}(a_1{\bf d}_1 + a_2{\bf d}_2 + \cdots + a_K{\bf d}_K) \\
		&= {\bf d}_k^\top {\bf A}_ka_k{\bf d}_k.
	\end{aligned}
\end{equation*}
So the coefficient can be obtained
\begin{equation}
	\label{eq-B-ak}
	a_k = \frac{{\bf d}_k^\top {\bf A}({\bf x}^* - {\bf x}_0)}{{\bf d}_k^\top {\bf A}_k{\bf d}_k}.
\end{equation}

Assume ${\bf x}_k$ is generated by the iterative scheme, then
\begin{equation*}
	\label{eq-B-xk}
	{\bf x}_{k-1} = {\bf x}_0 + \alpha_1{\bf d}_1 + \cdots \alpha_{k-1}{\bf d}_{k-1}.
\end{equation*}
Similarly, multiplying both sides of the formula by ${\bf d}_k^\top {\bf A}$, there is
\begin{equation}
	\label{eq-B-dkA(xk-x)}
	{\bf d}_k^\top {\bf A}({\bf x}_{k-1} - {\bf x}_0) = 0.
\end{equation}
Therefore,
\begin{equation}
	\label{eq-B-pkrk}
	\begin{aligned}
		{\bf d}_k^\top {\bf A} ({\bf x}^* - {\bf x}_0)
		&= {\bf d}_k^\top {\bf A} ({\bf x}^* - {\bf x}_{k-1}) \\
		& = {\bf d}_k^\top{\bf A} {\bf x}^* - {\bf A}{\bf x}_{k-1}) \\
		&= {\bf d}_k^\top({\bf b} - {\bf A}{\bf x}_{k-1}).
	\end{aligned}
\end{equation}

It has been mentioned in section \ref{sect2-1} that
\begin{equation*}
	\alpha_k = \frac{b_k - {\bf A}_k{\bf x}_{k-1} }{{\bf A}_k {\bf d}_k}.
\end{equation*}
Thus $b_k - {\bf A}_k{\bf x}_{k-1} = \alpha_k {\bf A}_k {\bf d}_k$ and substituting it into (\ref{eq-B-pkrk}), then
\begin{equation}
	\begin{aligned}
		{\bf d}_k^\top {\bf A} ({\bf x}^* - {\bf x}_0)
		&= {\bf d}_k^\top({\bf b} - {\bf A}{\bf x}_{k-1}) \\
		&= {\bf d}_k^\top\alpha_k {\bf A}{\bf d}_k \\
		&= {\bf d}_k^\top {\bf A} a_k{\bf d}_k.
	\end{aligned}
\end{equation}
That means $a_k = \alpha_k$, giving the result.

\section*{Appendix C.}
\label{Appendix C}
\setcounter{section}{3}
\setcounter{table}{0}
\setcounter{figure}{0}
\setcounter{equation}{0}
\renewcommand\theequation{C.\arabic{equation}}
\renewcommand{\thetable}{C.\arabic{table}}
\renewcommand{\thefigure}{C.\arabic{figure}}
This section gives a variant of the SOMA method and its convergence proof can refer to \cite{ref38}.

\begin{figure}[htbp]
	\centering
	\includegraphics[width=8.5cm]{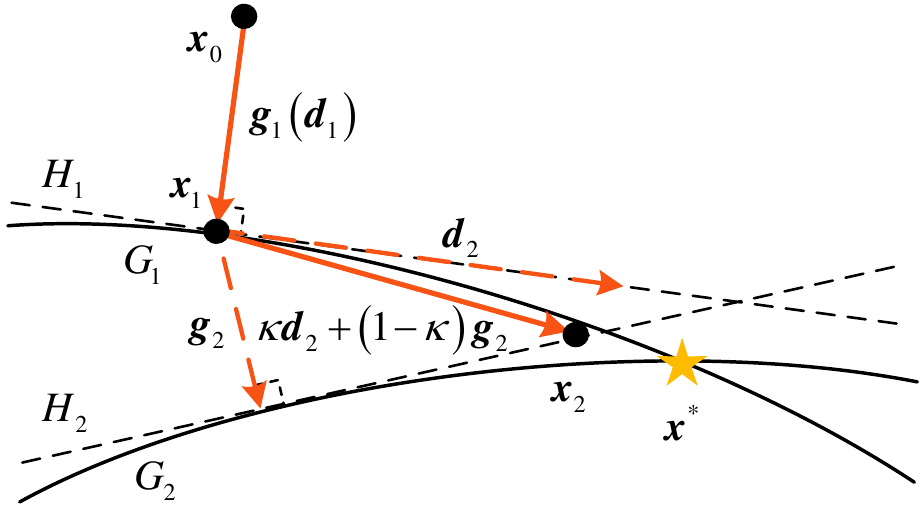} \hspace{-1cm}
	\vspace{-2.5mm}
	\caption{ The simplified geometric illustration of the variant method. The true solution ${\bf x}^*$ is represented by the yellow pentagram. $G_1$ and $G_2$ are two surfaces, and $H_1$ and $H_2$ are the corresponding tangent planes obtained by the first-order Taylor expansion.}
	\label{figC}
\end{figure}

The geometric illustration is shown in figure \ref{figC}. In the case of noise disturbance, the tangent planes can not approximate the surfaces. Thus, the change of the variant method is weighting the normal direction ${\bf g}_k$ and the orthogonal direction ${\bf d}_k$. The final search direction is $\kappa{\bf d}_k + (1-\kappa){\bf g}_k$, where $\kappa(\in [0,1])$ is the weighting coefficient.

By introducing weighting coefficient and relaxing the step size mentioned in section \ref{sect2-2}, the proposed method achieves a balance between the convergence speed and the solution accuracy, which can avoid the serious influence of noise amplification. In addition, the parameters of the proposed method have clear meanings, such as the weighting coefficient $\kappa$ and the step size relaxation factor $\beta$.

At the end of this section, the implementation of $\kappa$ is introduced in detail. For matched nonlinear equations without noise, the tangent planes can approximate the surfaces; in this case, $\kappa=1$ (or $\kappa\rightarrow 1$) and $\beta = 1$(or $\beta \rightarrow 1$). For mismatched nonlinear equations without noise, the constant terms are unknown; at this moment, the small step size can be used to search for solutions along ${\bf d}_k$, that is let $\kappa=1$ (or $\kappa\rightarrow 1$) and $\beta = 0.1$ (or $\beta \rightarrow 0$). For the case of noise disturbance, the corresponding linear equations are affected by the noise; at this time, that is $\kappa = 0.9$ or less and $\beta = 0.1$ (or $\beta \rightarrow 0$).

\section*{Appendix D.}
\label{Appendix D}
\setcounter{section}{4}
\setcounter{table}{0}
\setcounter{figure}{0}
\setcounter{equation}{0}
\renewcommand\theequation{D.\arabic{equation}}
\renewcommand{\thetable}{D.\arabic{table}}
\renewcommand{\thefigure}{D.\arabic{figure}}
The model of DSCT is mentioned in section \ref{sect1}. Sample at $\omega = 30,40,120,130$ and the sampling values of the mass attenuation coefficient (MAC) and the spectrum at corresponding energies are shown in table \ref{table-D-1}.

\begin{table}[]
	\footnotesize
	\centering
	\renewcommand\arraystretch{1.1}
	\setlength\tabcolsep{10pt}
	\caption{Some sampling values at specific energies.}
	\begin{tabular}{cccc}
		\hline
		Energy & Bone MAC & Water MAC & Spectrum sampling value \\ \hline
		30 & 0.2812 & 0.0395 & 0.0002 \\
		40 & 0.1342 & 0.0281 & 0.0009 \\
		120 & 0.0328 & 0.0159 & 0.0056 \\
		130 & 0.0314 & 0.0154 & 0.0029 \\ \hline
	\end{tabular}
	\label{table-D-1}
\end{table}

The first spectrum consists of 30 keV and 40 keV and the second spectrum consists of 120 keV and 130 keV. Then, a simplified DSCT example can be express as
\begin{equation}
	\label{eq-D-simDSCT}
	\left\{
	\begin{aligned}
		p_{1} &= -\ln(0.0002e^{-0.2812q_1-0.0395q_2} + 0.0009e^{-0.1342q_1-0.0281q_2}), \\
		p_{2} &= -\ln(0.0056e^{-0.0328q_1-0.0159q_2} + 0.0029e^{-0.0314q_1-0.0154q_2}).
	\end{aligned}
	\right.
\end{equation}

\begin{figure}[htbp]
	\centering
	\includegraphics[height=7.5cm]{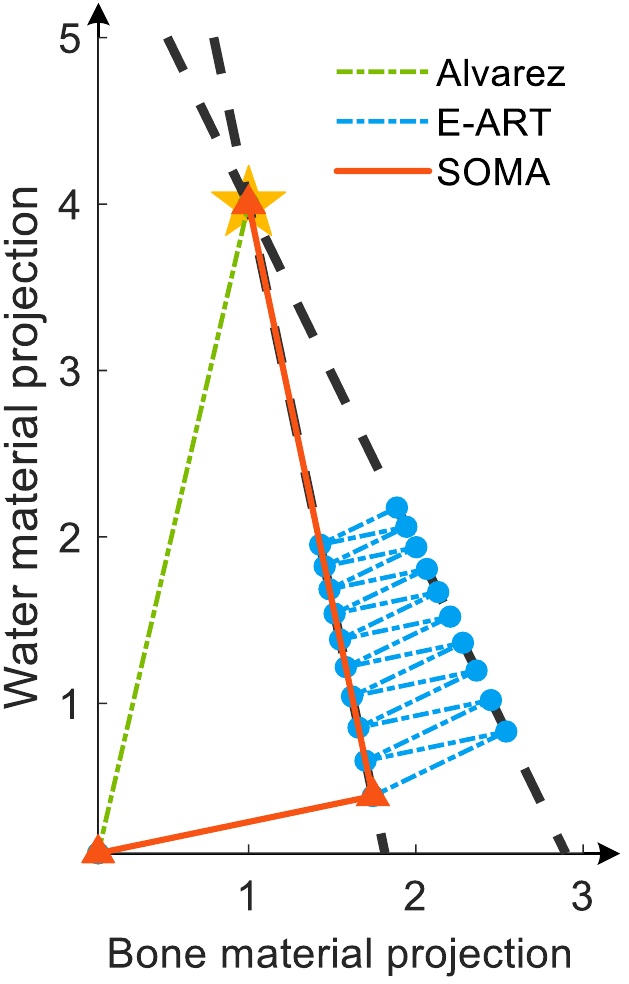} \hspace{-1cm}
	\vspace{-2.5mm}
	\caption{ Iterative paths of some methods for solving the simplified DSCT example. The true solution $(1,4)$ is represented by the yellow pentagram. The black dashed curves represent the nonlinear equations (though they look like lines).}
	\label{figD}
\end{figure}

Suppose the true solutions are $q1 = 1$ and $q2 = 4$. The iterative paths of the Alvarez method, the E-ART method and the proposed method for solving the nonlinear equations (\ref{eq-D-simDSCT}) are shown in figure \ref{figD}. Because the iterative paths of the NCPD method are similar to those of the E-ART method, no iterative paths of it are plotted to show clearly. 

It can be observed that the Alvarez method and the proposed method have fast convergence speeds. But the Alvarez method only can be applied to matched nonlinear equations and is not robust to noise. On the other hand, the proposed method only needs two steps to get the true solution, while the E-ART method still can not reach the neighborhood of the true solution after many iterations.

\section*{Reference}
\label{sect5}

\end{document}